\documentclass[final,3p,times]{elsarticle}
\usepackage{CJK}
\usepackage{multirow}
\usepackage{booktabs}
\usepackage{amssymb}
\usepackage{comment}
\usepackage{mathrsfs}        
\usepackage{amsthm}          
\usepackage{amsmath}         
\usepackage{epstopdf}
\biboptions{numbers,sort&compress}
\usepackage[colorlinks, citecolor=blue]{hyperref}
\hypersetup{colorlinks=false}                       
\usepackage{xcolor}  
\definecolor{b}{rgb}{0,0,1}
\usepackage{caption2}                   
\theoremstyle{plain}

\theoremstyle{definition}
\newtheorem{definition}{Definition}
\newtheorem{example}{Example}
\newtheorem{remark}{Remark}

\journal{Applied Mathematics and Computation}

\begin{document}

\begin{frontmatter}



\title{Infinite series representation of fractional calculus: theory and applications}

\author[ustc]{Yiheng Wei}
\author[ucm]{YangQuan Chen}
\author[ude]{Qing Gao}
\author[ustc]{Yong Wang\corref{cor1}}
\ead{yongwang@ustc.edu.cn}
\cortext[cor1]{Corresponding author }
\address[ustc]{Department of Automation, University of Science and Technology of China, Hefei, 230026, China}
\address[ucm]{School of Engineering, University of California, Merced,  Merced, 95343, USA}
\address[ude]{The Institute for Automatic Control and Complex Systems, University of Duisburg-Essen, Duisburg 47057, Germany}

\begin{abstract}
This paper focuses on the equivalent expression of fractional integrals/derivatives with an infinite series. A universal framework for fractional Taylor series is developed by expanding an analytic function at the initial instant or the current time. The framework takes into account of the Riemann--Liouville definition, the Caputo definition, the constant order and the variable order. On this basis, some properties of fractional calculus are confirmed conveniently. An intuitive numerical approximation scheme via truncation is proposed subsequently. Finally, several illustrative examples are presented to validate the effectiveness and practicability of the obtained results.
\end{abstract}


\end{frontmatter}



\section{Introduction}\label{Section 1}
Taylor series has intensively developed since its introduction in 300 years ago and is nowadays a mature research field. As a powerful tool, Taylor series plays an essential role in analytical analysis and numerical calculation of a function. Interestingly, the classical Taylor series has been tied to another 300-year history tool, i.e., fractional calculus. This combination produces many promising and potential applications \cite{Osler:1971JMA,Chen:2017AMC,Jaradat:2018CSF,Fernandez:2019CNSNS}.

The original idea on fractional generalized Taylor series can be dated back to 1847, when Riemann formally used a series structure to formulate an analytic function \cite{Riemann:1953Book}. This series was not proven and the related manuscript probably never intended for publishing. The proof of the validity of such an expansion for certain classes of functions was given by Hardy \cite{Hardy:1945JLMS}, both for finite and infinite initial instant. Following these pioneer work, Trujillo et al. discussed the related mean value theorem problem \cite{Trujillo:1999JMAA} and Odibat et al. extended the result from the Riemann--Liouville case to the Caputo case \cite{Odibat:2007AMC}. Since then, a large volume of papers have been published sequentially on this topic, to mention a few, see \cite{Odibat:2008AMC,Tremblay:2013ITSF,Krishnasamy:2017JAS,Benjemaa:2018AMC,Fernandez:2018ADE,Cheng:2019JMS}. However, these results mainly focused on how to describe a function or a fractional derivative with a series comprised of fractional derivatives.

Luckily, Oldham and Spanier suggested a different perspective to express fractional calculus \cite{Oldham:1974Book}, including Gr\"{u}nwald--Letnikov and Riemann--Liouville definitions. In this expression, a power series involving integer derivatives of the discussed analytic function was constructed. Apart from this series, Samko et al. established a new series to describe fractional derivative \cite{Samko:1993Book}. The two series provide a particular insight into an interesting research area and opens the reader's mind to a world yet completely unexplored. Afterwards, the series representation of fractional calculus with fixed memory length \cite{Wei:2017FCAAb} and discrete time \cite{Wei:2019CNSNS,Wei:2019arXiv} were investigated. Besides, some preprint papers are found on the theory and application of fractional series representation \cite{Gladkina:2018arXiv,Sousa:2018arXiv,Shchedrin:2018arXiv}.

Despite of some pioneer work reported, as the authors know, it has not yet attracted enough attention in both engineering and non-engineering disciplines. i) A compact and practical series representation is desired for the Caputo fractional derivative. ii) Many properties are desired to be explore with the developed series. iii) Compared with the constant order case, the variable case exhibits more evidence on the natural property of fractional calculus, while related research is absent. Bearing these ideas in mind, this paper aims at deriving infinite series representations for fractional integrals/derivatives and addressing their underlying properties comprehensively.

The remainder of the paper is divided in the following manner. Section \ref{Section 2} provides some preliminaries of fractional calculus. The theory and applications of series representation for fractional integrals/derivaters are investigated in Section \ref{Section 3}. In Section \ref{Section 4}, three numerical examples are presented to illustrate the validity of the developed method. This paper is concluded in Section \ref{Section 5}.

\section{Preliminaries}\label{Section 2}
Fractional calculus could be treated as a natural generalization of classical integer calculus. There have been a number of different definitions for the fractional derivative, among which the Riemann--Liouville fractional derivative and the Caputo fractional derivative are most widely used.

\begin{definition} \label{Definition 1}
The $\alpha$-th Riemann--Liouville fractional integral of a function $f(t)$ is defined as
\begin{equation}\label{Eq1}
{\textstyle
{}_a^{\rm R}{{\mathscr I}_{t}^\alpha}{f\left( t \right)}  \triangleq \frac{1}{{\Gamma \left( \alpha  \right)}}\int_{{a}}^t {{{\left( {t - \tau } \right)}^{\alpha -1}}{ f\left( \tau  \right)}{\rm{d}}\tau },}
\end{equation}
where $\alpha>0$ is the integral order, $a$ is the initial instant and $\Gamma \left( x \right) = \int_0^{ + \infty } {{{\rm{e}}^{ - t}}{t^{x - 1}}{\rm{d}}t} $ is the Gamma function.
\end{definition}

In the light of the fractional integral in (\ref{Eq1}), the following two different fractional derivatives could be established respectively \cite{Oldham:1974Book}. Let $n\in\mathbb{N}_+$ and $\alpha\in(n-1,n)$.

\begin{definition} \label{Definition 2}
The $\alpha$-th Riemann--Liouville fractional derivative of a function $f(t)$ is defined as ${}^{\rm R}_a{\mathscr D}_{t}^\alpha f\left( t \right) \triangleq \frac{{\rm{d}}^n}{{{\rm{d}}t^n}} {_a^{\rm R}{\mathscr I}_{t}^{n-\alpha}} f\left( t \right)$.
\end{definition}

\begin{definition} \label{Definition 3}
The $\alpha$-th Caputo fractional derivative of a function $f(t)$ is defined as ${}^{\rm C}_a{\mathscr D}_{t}^\alpha f\left( t \right) \triangleq  {}_a^{\rm R}{\mathscr I}_{t}^{n-\alpha}\frac{{\rm{d}}^n}{{{\rm{d}}t^n}} f\left( t \right)$.
\end{definition}

The singularity issue of fractional calculus can be directly concluded from (\ref{Eq1}), since there will be zero denominator in (\ref{Eq1}) as $\tau\to t$ and $0<\alpha<1$. Similarly, for the derived derivatives, $0<n-\alpha<1$ also results in singularity. One can also read from the above definitions that a fractional derivative is actually a special fractional integral that depends on all the historical data. This nonlocality is also called the \textit{global correlation} or \textit{history dependence} \cite{Sun:2017PSMA}, which signifies the essential difference between fractional calculus and classical integer calculus.

\section{Main Results}\label{Section 3}
This section proposes two series representations for fractional derivative and fractional integral, based on which several calculation properties of fractional calculus are deduced.

\subsection{A universal framework}
If $f(\cdot)$ is an analytic function in an interval $(a,b)$, then it can be represented as a convergent series $f(t)=\sum\nolimits_{k = 0}^{ + \infty } {\frac{{{f^{\left( k \right)}}\left( a \right)}}{{k!}}{{\left( {t - a} \right)}^k}} $. With this assumption, the $\alpha$-th Caputo fractional derivative of $f(t)$ can be calculated as
 \begin{equation}\label{Eq2}
{\textstyle
\begin{array}{rl}
{}^{\rm C}_a{\mathscr D}_{t}^\alpha f\left( t \right) =&\hspace{-6pt}{}_a^{\rm R}{\mathscr I}_{t}^{n-\alpha}\frac{{\rm{d}}^n}{{{\rm{d}}t^n}} f\left( t \right)\\
=&\hspace{-6pt}{}_a^{\rm R}{\mathscr I}_{t}^{n-\alpha}\frac{{\rm{d}}^n}{{{\rm{d}}t^n}} {\sum\nolimits_{k = 0}^{ + \infty } {\frac{{{f^{\left( k \right)}}\left( a \right)}}{{k!}}{{\left( {t - a} \right)}^k}}}\\
=&\hspace{-6pt}{}_a^{\rm R}{\mathscr I}_{t}^{n-\alpha} {\sum\nolimits_{k = n}^{ + \infty } {\frac{{{f^{\left( k \right)}}\left( a \right)}}{{\Gamma\left(k-n+1\right)}}{{\left( {t - a} \right)}^{k-n}}}}\\
=&\hspace{-6pt}{\sum\nolimits_{k = n}^{ + \infty } {\frac{{{f^{\left( k \right)}}\left( a \right)}}{{\Gamma\left(k-n+1\right)}}{}_a^{\rm R}{\mathscr I}_{t}^{n-\alpha}{{\left( {t - a} \right)}^{k-n}}}}\\
=&\hspace{-6pt} \sum\nolimits_{k = n}^{+\infty}  {\frac{{{f^{\left( k \right)}}\left( a \right)}}{{\Gamma \left( { k- \alpha + 1 } \right)}}} {\left( {t - a} \right)^{k - \alpha }},
\end{array}}
\end{equation}
where $n-1<\alpha<n$ and $n\in \mathbb{N}_+$.

In (\ref{Eq2}), the Taylor series is expanded at the initial instant $t=a$. The corresponding case expanded at the current time comes to mind intuitively. If $f(\cdot)$ can be expressed a convergent series $f(\tau)=\sum\nolimits_{k = 0}^{ + \infty } {\frac{{{f^{\left( k \right)}}\left( t \right)}}{{k!}}{{\left( {\tau - t} \right)}^k}} $, $\tau\in(a,t)$, then another useful series representation can be achieved regarding the following deduction
\begin{equation}\label{Eq3}
\begin{array}{rl}
{}_a^{\rm{C}}{\mathscr D}_t^\alpha f\left( t \right)=&\hspace{-6pt}\frac{1}{{\Gamma \left( {n - \alpha } \right)}}\int_a^t {{{\left( {t - \tau } \right)}^{n - \alpha  - 1}}{f^{\left( n \right)}}\left( \tau  \right){\rm{d}}\tau } \\
=&\hspace{-6pt}\frac{1}{{\Gamma \left( {n - \alpha } \right)}}\int_a^t {\sum\nolimits_{k = 0}^{ + \infty } {{{\left( { - 1} \right)}^k}\frac{{{f^{\left( {k + n} \right)}}\left( t \right)}}{{\Gamma \left( {k + 1} \right)}}} {{\left( {t - \tau } \right)}^{n + k - \alpha  - 1}}{\rm{d}}\tau } \\
=&\hspace{-6pt}\frac{1}{{\Gamma \left( {n - \alpha } \right)}}\sum\nolimits_{k = 0}^{ + \infty } {{{\left( { - 1} \right)}^k}\frac{{{f^{\left( {k + n} \right)}}\left( t \right)}}{{\Gamma \left( {k + 1} \right)}}} \int_a^t {{{\left( {t - \tau } \right)}^{n + k - \alpha  - 1}}{\rm{d}}\tau } \\
=&\hspace{-6pt}\frac{1}{{\Gamma \left( {n - \alpha } \right)}}\sum\nolimits_{k = 0}^{ + \infty } {{{\left( { - 1} \right)}^k}\frac{{{f^{\left( {k + n} \right)}}\left( t \right)}}{{\Gamma \left( {k + 1} \right)}}} \frac{{{{\left( {t - a} \right)}^{n + k - \alpha }}}}{{n + k - \alpha }}\\
=&\hspace{-6pt} \sum\nolimits_{k = n}^{ + \infty } {\frac{{{{\left( { - 1} \right)}^k}{f^{\left( k \right)}}\left( t \right)}}{{\Gamma \left( {n - \alpha } \right)\Gamma \left( {k - n - 1} \right)\left( {k - \alpha } \right)}}{{\left( {t - a} \right)}^{k - \alpha }}} \\
=&\hspace{-6pt} \sum\nolimits_{k = n}^{ + \infty } {\left( {\begin{smallmatrix}
{\alpha  - n}\\
{k - n}
\end{smallmatrix}} \right)\frac{{{f^{\left( k \right)}}\left( t \right)}}{{\Gamma \left( {k - \alpha+ 1} \right)}}{{\left( {t - a} \right)}^{k - \alpha }}},
\end{array}
\end{equation}
where $n - 1 < \alpha  < n$, $n\in\mathbb{N}_+$, the reflection formula $\Gamma \left( x \right)\Gamma \left( {1 - x} \right) = \frac{\pi }{{\sin \left( {x\pi } \right)}}$ is used and the binomial coefficient is defined as $\left( {\begin{smallmatrix}
p \\
q
\end{smallmatrix}} \right) = \frac{{\Gamma \left( {p  + 1} \right)}}{{\Gamma \left( {q  + 1} \right)\Gamma \left( {p  - q  + 1} \right)}}$.

Notably, from another perspective, the first class in (\ref{Eq2}) is derived by expanding $f(\cdot)$ before substituting into the specific definition and the second case in (\ref{Eq3}) is established by expanding $f(\cdot)$ after substituting into the specific definition. Note that with the adoption of $f^{(k)}(t)$, the second class is not a power series like the first class.

Similar to the newly developed series representation in Caputo case, some existing results emerged in pioneering work (see (3.53), (5.2.4) and (5.2.9) of \cite{Oldham:1974Book}; (15.4) of \cite{Samko:1993Book}; (2.215) and (2.223) of \cite{Podlubny:1999Book})
\begin{equation}\label{Eq4}
{\textstyle{}_a^{\rm R}{\mathscr I}_t^\alpha f\left( t \right) = \sum\nolimits_{k = 0}^{ + \infty } {\frac{{{f^{\left( k \right)}}\left( a \right)}}{{\Gamma \left( {k + \alpha+ 1 } \right)}}{{\left( {t - a} \right)}^{k + \alpha }}},}
\end{equation}
\begin{equation}\label{Eq5}
{\textstyle{}^{\rm R}_a{\mathscr D}_{t}^\alpha f\left( t \right) =\sum\nolimits_{k = 0}^{+\infty}  {\frac{{{f^{\left( k \right)}}\left( a \right)}}{{\Gamma \left( { k- \alpha + 1 } \right)}}} {\left( {t - a} \right)^{k - \alpha }},}
\end{equation}
\begin{equation}\label{Eq6}
{\textstyle{}_a^{\rm R}{\mathscr I}_t^\alpha f\left( t \right) =\sum\nolimits_{k = 0}^{ + \infty } {\left( {\begin{smallmatrix}
{ - \alpha }\\
k
\end{smallmatrix}} \right)\frac{{{f^{\left( k \right)}}\left( t \right)}}{{\Gamma \left( {k + \alpha+ 1} \right)}}{{\left( {t - a} \right)}^{k + \alpha }}},}
\end{equation}
\begin{equation}\label{Eq7}
{\textstyle{}_a^{\rm{R}}{\mathscr D}_t^\alpha f\left( t \right)= \sum\nolimits_{i = 0}^{ + \infty } {\left( {\begin{smallmatrix}
\alpha \\
i
\end{smallmatrix}} \right)\frac{{{f^{\left( i \right)}}\left( t \right)}}{{\Gamma \left( {i - \alpha+ 1 } \right)}}{{\left( {t - a} \right)}^{i - \alpha }}} ,}
\end{equation}
where similar assumption on $f(t)$ is made, $\alpha>0$ and $t\in[a,b]$.

Interestingly, it can be observed that if we extend the range of $\alpha$ to $(-\infty,+\infty)$, (\ref{Eq4}) and (\ref{Eq5}) share a same form. Likewise, (\ref{Eq6}) and (\ref{Eq7}) can also share a unified form with $\alpha\in(-\infty,+\infty)$. Notably, such a fact can also be verified through the following deduction ($n-1<\alpha<n\in\mathbb{N}_+$)
\begin{equation}\label{Eq8}
\begin{array}{rl}
{}_a^{\rm{R}}{\mathscr D}_t^{-\alpha} f\left( t \right) =&\hspace{-6pt} \frac{{{{\rm{d}}^n}}}{{{\rm{d}}{t^n}}}\frac{1}{{\Gamma \left( {n + \alpha } \right)}}\int_a^t {{{\left( {t - \tau } \right)}^{n + \alpha  - 1}}f\left( \tau  \right){\rm{d}}\tau } \\
=&\hspace{-6pt} \frac{1}{{\Gamma \left( {  \alpha } \right)}}\int_a^t {{{\left( {t - \tau } \right)}^{  \alpha  - 1}}f\left( \tau  \right){\rm{d}}\tau } \\
 =&\hspace{-6pt} {}_a^{\rm R}{\mathscr I}_t^{  \alpha }f\left( t \right),
\end{array}
\end{equation}
where the formula $\frac{{\rm{d}}}{{{\rm{d}}t}}\int_a^t {f\left( {t,\tau } \right){\rm{d}}\tau}  = \int_a^t {\frac{\partial }{{\partial t}}f\left( {t,\tau } \right){\rm{d}}\tau}  + f\left( {t,t} \right)$ is used repeatedly.

Comparing (\ref{Eq2}) with (\ref{Eq5}), a well-known formula ($n-1<\alpha<n\in\mathbb{N}_+$) follows
\begin{equation}\label{Eq9}
{\textstyle
{}^{\rm R}_a{\mathscr D}_{t}^\alpha f\left( t \right) = {}^{\rm C}_a{\mathscr D}_{t}^\alpha f\left( t \right)+\sum\nolimits_{k = 0}^{n-1}  {\frac{{{f^{\left( k \right)}}\left( a \right)}}{{\Gamma \left( { k- \alpha + 1 } \right)}}} {\left( {t - a} \right)^{k - \alpha }}.}
\end{equation}
With the representations in (\ref{Eq3}) with (\ref{Eq7}), the relationship between ${}_a^{\rm{R}}{\mathscr D}_t^\alpha f\left( t \right)$ and ${}_a^{\rm{C}}{\mathscr D}_t^\alpha f\left( t \right)$ cannot be expressed as a brief formula like (\ref{Eq9}).

To improve the practicality, we reduce the requirement on $f(\cdot)$. If the function $f\in{\cal C}^{N+1}[a,b]$ can be expressed by a Taylor formula expanded at the initial instant, i.e., $f\left( t \right) = \sum\nolimits_{k = 0}^{ N} {\frac{{{f^{\left( k \right)}}\left( a \right)}}{k!}{{\left( {t - a} \right)}^{k }}}  +\frac{1}{\Gamma(N+1)}\int_a^t {{{\left( {t - \tau } \right)}^{N}}f^{(N+1)}\left( \tau  \right){\rm{d}}\tau }$, $N\in\mathbb{N}$, then a new equivalent form of ${}_a^{\rm R}{\mathscr I}_t^\alpha f\left( t \right)$ can be obtained
\begin{equation}\label{Eq10}
\begin{array}{rl}
{}_a^{\rm R}{\mathscr I}_t^\alpha f\left( t \right) =&\hspace{-6pt} {}_a^{\rm R}{\mathscr I}_t^\alpha\sum\nolimits_{k = 0}^N {\frac{{{f^{\left( k \right)}}\left( a \right)}}{{k!}} {{\left( {t - a} \right)}^k}}  + {}_a^{\rm R}{\mathscr I}_t^\alpha\frac{1}{{\Gamma(N + 1)}} \int_a^t {{{\left( {t - \tau } \right)}^N}{f^{(N + 1)}}\left( \tau  \right){\rm{d}}\tau } \\
 =&\hspace{-6pt} \sum\nolimits_{k = 0}^N {\frac{{{f^{\left( k \right)}}\left( a \right)}}{{\Gamma \left( {k + \alpha  + 1} \right)}}{{\left( {t - a} \right)}^{k + \alpha }}}  + \frac{1}{{\Gamma \left( \alpha  \right)}}\int_a^t {{{\left( {t - x } \right)}^{\alpha  - 1}}} \frac{1}{{\Gamma (N + 1)}}\int_a^x  {{{\left( {x  - \tau} \right)}^N}{f^{(N + 1)}}\left( \tau \right){\rm{d}}\tau{\rm{d}}x }.
\end{array}
\end{equation}

The second term in the right hand of (\ref{Eq10}) can be further simplified as
\begin{equation}\label{Eq11}
\begin{array}{l}
\frac{1}{{\Gamma \left( \alpha  \right)}}\int_a^t {{{\left( {t - x} \right)}^{\alpha  - 1}}} \frac{1}{{\Gamma (N + 1)}}\int_a^x {{{\left( {x - \tau } \right)}^N}{f^{(N + 1)}}\left( \tau  \right){\rm{d}}\tau {\rm{d}}x} \\
 = \frac{1}{{\Gamma \left( \alpha  \right)\Gamma (N + 1)}}\int_a^t {\int_a^x {{{\left( {t - x} \right)}^{\alpha  - 1}}{{\left( {x - \tau } \right)}^N}{f^{(N + 1)}}\left( \tau  \right){\rm{d}}\tau {\rm{d}}x} } \\
 = \frac{1}{{\Gamma \left( \alpha  \right)\Gamma (N + 1)}}\int_a^t {\int_\tau ^t {{{\left( {t - x} \right)}^{\alpha  - 1}}{{\left( {x - \tau } \right)}^N}{f^{(N + 1)}}\left( \tau  \right){\rm{d}}x{\rm{d}}\tau } } \\
 = \frac{1}{{\Gamma \left( \alpha  \right)\Gamma (N + 1)}}\int_a^t {\int_\tau ^t {{{\left( {t - x} \right)}^{\alpha  - 1}}{{\left( {x - \tau } \right)}^N}{\rm{d}}x} {f^{(N + 1)}}\left( \tau  \right)} {\rm{d}}\tau \\
 = \frac{1}{{\Gamma \left( \alpha  \right)\Gamma (N + 1)}}\int_a^t {\int_0^1 {{{\left( {t - \tau } \right)}^{N + \alpha }}{{\left( {1 - \xi } \right)}^{\alpha  - 1}}{\xi ^N}{\rm{d}}\xi } {f^{(N + 1)}}\left( \tau  \right)} {\rm{d}}\tau \\
 = \frac{1}{{\Gamma (N + \alpha  + 1)}}\int_a^t {{{\left( {t - \tau } \right)}^{N + \alpha }}{f^{(N + 1)}}\left( \tau  \right)} {\rm{d}}\tau.
\end{array}
\end{equation}
Substituting (\ref{Eq11}) into (\ref{Eq10}), a beautiful result can be obtained
\begin{equation}\label{Eq12}
{\textstyle {}_a^{\rm R}{\mathscr I}_t^\alpha f\left( t \right) = \sum\nolimits_{k = 0}^{ N} {\frac{{{f^{\left( k \right)}}\left( a \right)}}{{\Gamma \left( {k+ \alpha + 1 } \right)}}{{\left( {t - a} \right)}^{k + \alpha }}} +\frac{1}{{\Gamma \left( N+\alpha+1  \right)}}\int_a^t {{{\left( {t - \tau } \right)}^{N+\alpha}}f^{(N+1)}\left( \tau  \right){\rm{d}}\tau },}
\end{equation}
where $\alpha>0$ and $t\in[a,b]$.

Considering $\alpha \in(n-1,n)$ and $n\in\mathbb{N}_+$, by using Definition \ref{Definition 2}, one has
\begin{equation}\label{Eq13}
{\textstyle \begin{array}{rl}
{}_a^{\rm R}{\mathscr D}_t^\alpha f\left( t \right) =&\hspace{-6pt}\frac{{{{\rm{d}}^n}}}{{{\rm{d}}{t^n}}}\sum\nolimits_{k = 0}^N {\frac{{{f^{\left( k \right)}}\left( a \right)}}{{\Gamma \left( {k + n - \alpha  + 1} \right)}}{{\left( {t - a} \right)}^{k + n - \alpha }}}  + \frac{{{{\rm{d}}^n}}}{{{\rm{d}}{t^n}}}\frac{1}{{\Gamma \left( {N + n - \alpha  + 1} \right)}}\int_a^t {{{\left( {t - \tau } \right)}^{N + n - \alpha }}{f^{(N + 1)}}\left( \tau  \right){\rm{d}}\tau } \\
 =&\hspace{-6pt} \sum\nolimits_{k = 0}^N {\frac{{{f^{\left( k \right)}}\left( a \right)}}{{\Gamma \left( {k - \alpha  + 1} \right)}}{{\left( {t - a} \right)}^{k - \alpha }}}  + \frac{{{{\rm{d}}^{n - 1}}}}{{{\rm{d}}{t^{n - 1}}}}\frac{1}{{\Gamma \left( {N + n - \alpha } \right)}}\int_a^t {{{\left( {t - \tau } \right)}^{N + n - \alpha  - 1}}{f^{(N + 1)}}\left( \tau  \right){\rm{d}}\tau } \\
 &\hspace{-6pt}+ \frac{{{{\rm{d}}^{n - 1}}}}{{{\rm{d}}{t^{n - 1}}}}\frac{{{{\left( {t - t} \right)}^{N + n - \alpha }}{f^{(N + 1)}}\left( t \right)}}{{\Gamma \left( {N + n - \alpha  + 1} \right)}}\\
 =&\hspace{-6pt} \sum\nolimits_{k = 0}^N {\frac{{{f^{\left( k \right)}}\left( a \right)}}{{\Gamma \left( {k - \alpha  + 1} \right)}}{{\left( {t - a} \right)}^{k - \alpha }}}  + \frac{1}{{\Gamma \left( {N - \alpha  + 1} \right)}}\int_a^t {{{\left( {t - \tau } \right)}^{N - \alpha }}{f^{(N + 1)}}\left( \tau  \right){\rm{d}}\tau },
\end{array}}
\end{equation}
where the differential formula for integral upper limit function is adopted repeatedly.

When $N\in\mathbb{N}$ and $N\ge n$, the Caputo case can be derived
\begin{equation}\label{Eq14}
{\textstyle
\begin{array}{rl}
{}_a^{\rm C}{\mathscr D}_t^\alpha f\left( t \right) =&\hspace{-6pt}
{}_a^{\rm{R}}{\mathscr I}_t^{n - \alpha }\frac{{{{\rm{d}}^n}}}{{{\rm{d}}{t^n}}}\sum\nolimits_{k = 0}^N {\frac{{{f^{\left( k \right)}}\left( a \right)}}{{k!}}{{\left( {t - a} \right)}^k}}  + {}_a^{\rm{R}}{\mathscr I}_t^{n - \alpha }\frac{{{{\rm{d}}^n}}}{{{\rm{d}}{t^n}}}\frac{1}{{\Gamma (N + 1)}}\int_a^t {{{\left( {t - \tau } \right)}^N}{f^{(N + 1)}}\left( \tau  \right){\rm{d}}\tau } \\
 =&\hspace{-6pt} {}_a^{\rm{R}}{\mathscr I}_t^{n - \alpha }\sum\nolimits_{k = n}^N {\frac{{{f^{\left( k \right)}}\left( a \right)}}{{\Gamma \left( {k - n + 1} \right)}}{{\left( {t - a} \right)}^{k - n}}}  + {}_a^{\rm{R}}{\mathscr I}_t^{n - \alpha }\frac{1}{{\Gamma (N - n + 1)}}\int_a^t {{{\left( {t - \tau } \right)}^{N - n}}{f^{(N + 1)}}\left( \tau  \right){\rm{d}}\tau } \\
 =&\hspace{-6pt}\sum\nolimits_{k = n}^{ N} {\frac{{{f^{\left( k \right)}}\left( a \right)}}{{\Gamma \left( {k - \alpha+ 1 } \right)}}{{\left( {t - a} \right)}^{k - \alpha }}}+\frac{1}{{\Gamma \left( N-\alpha+1  \right)}}\int_a^t {{{\left( {t - \tau } \right)}^{N-\alpha}}f^{(N+1)}\left( \tau  \right){\rm{d}}\tau }.
\end{array}}
\end{equation}
Note that the remainder in (\ref{Eq14}) is equal to that in (\ref{Eq13}), which confirms the relationship in (\ref{Eq9}).

Likewise, if the function $f\in{\cal C}^{N+1}[a,b]$ can be expressed by a Taylor formula expanded at the current time $\tau =t$, i.e., $f\left( \tau \right) = \sum\nolimits_{k = 0}^N {\frac{{{f^{\left( k \right)}}\left( t \right)}}{{\Gamma \left( {k + 1} \right)}}} {\left( {\tau - t} \right)^k} + \frac{1}{{\Gamma \left( {N + 1} \right)}}\int_a^{\tau} {{{\left( {\tau-x } \right)}^N}{f^{\left( {N + 1} \right)}}\left( x  \right){\rm{d}}x }$, $\tau\in(a,t)$, one has
\begin{equation}\label{Eq15}
{\textstyle {}_a^{\rm R}{\mathscr I}_t^\alpha f\left( t \right) = \sum\nolimits_{k = 0}^{ N } {\left( {\begin{smallmatrix}
{ - \alpha }\\
k
\end{smallmatrix}} \right)\frac{{{f^{\left( k \right)}}\left( t \right)}}{{\Gamma \left( {k+ \alpha+ 1} \right)}}{{\left( {t - a} \right)}^{k + \alpha }}}  +\frac{1}{{\Gamma \left( N+\alpha+1  \right)}}\int_a^t {{{\left( {t - \tau } \right)}^{N+\alpha}}f^{(N+1)}\left( \tau  \right){\rm{d}}\tau },}
\end{equation}
where $\alpha>0$ and $N\in\mathbb{N}$.

When $\alpha\in(n-1,n)$, $n\in\mathbb{N}_+$, $N\in\mathbb{N}$ and $N\ge n$, the following equations follow from the relationship between fractional derivative and fractional integral immediately
\begin{equation}\label{Eq16}
{\textstyle{}_a^{\rm R}{\mathscr D}_t^\alpha f\left( t \right) = \sum\nolimits_{k = 0}^{ N } {\left( {\begin{smallmatrix}
\alpha \\
k
\end{smallmatrix}} \right)\frac{{{f^{\left( k \right)}}\left( t \right)}}{{\Gamma \left( {k- \alpha + 1 } \right)}}{{\left( {t - a} \right)}^{k - \alpha }}}   +\frac{1}{{\Gamma \left( N-\alpha+1  \right)}}\int_a^t {{{\left( {t - \tau } \right)}^{N-\alpha}}f^{(N+1)}\left( \tau  \right){\rm{d}}\tau },}
\end{equation}
\begin{equation}\label{Eq17}
{\textstyle{}_a^{\rm C}{\mathscr D}_t^\alpha f\left( t \right) = \sum\nolimits_{k = n}^{N } {\left( {\begin{smallmatrix}
{\alpha  - n}\\
{k - n}
\end{smallmatrix}} \right)\frac{{{f^{\left( k \right)}}\left( t \right)}}{{\Gamma \left( {k- \alpha + 1} \right)}}{{\left( {t - a} \right)}^{k - \alpha }}} +\frac{1}{{\Gamma \left( N+n-\alpha+1  \right)}}\int_a^t {{{\left( {t - \tau } \right)}^{N+n-\alpha}}f^{(N+n+1)}\left( \tau  \right){\rm{d}}\tau }.}
\end{equation}

\begin{remark}\label{Remark 1}
By now, the universal representation framework has been developed for fractional integral and fractional derivative, which includes the Taylor series representation and the Taylor formula representation. Compared with the Taylor series based formulas (\ref{Eq2})-(\ref{Eq7}), lower requirements are needed for $f(\cdot)$ in the Taylor formula based ones (\ref{Eq12})-(\ref{Eq17}). As a result, the latter ones provide more practicability. On the other side, the new representation can be divided into two kinds: the one expanded at the initial instant and the one expanded at the current time.
\end{remark}

In order to obtain a comprehensive understanding of the elaborated framework, several applications are presented in the remainder of this section. To be specifically, subsections 3.2-3.4 analyzes a situation where the first class of series is more convenient to use. Subsection 3.5 shows that we may need different series representations in different analysis stages of the same problem. In subsection 3.6, we provide a case when both the two kinds series are useful. The results on constant order $\alpha$ are extended to the variable order case in subsection 3.7. Consider an analytic function $f(\cdot)$ in an interval $(a,b)$ unless otherwise specified.

\subsection{Associative law}
Assume $f(\cdot)$ can be expressed as a Taylor series at the initial instant, i.e., $f(t)=\sum\nolimits_{k = 0}^{ + \infty } {\frac{{{f^{\left( k \right)}}\left( a \right)}}{{k!}}{{\left( {t - a} \right)}^k}} $, $t\in[a,b]$. The associate laws of fractional derivative in two cases are formulated as follows.

\begin{itemize}
  \item Case 1: with fractional derivative.
\end{itemize}

Let ${\alpha _1},{\alpha _2}>0$ and $\alpha_1,\alpha_2\notin \mathbb{N}$, one has
\begin{eqnarray}\label{Eq18}
\begin{array}{rl}
{}_a^{\rm{R}}{\mathscr D}_t^{{\alpha _1}}{}_a^{\rm{R}}{\mathscr D}_t^{{\alpha _2}}f\left( t \right) = &\hspace{-6pt} {}_a^{\rm{R}}{\mathscr D}_t^{{\alpha _1}}\sum\nolimits_{k = 0}^{ + \infty } {\frac{{{f^{\left( k \right)}}\left( a \right)}}{{\Gamma \left( {k - {\alpha _2}+ 1} \right)}}{{\left( {t - a} \right)}^{k - {\alpha _2}}}} \\
 =&\hspace{-6pt} \sum\nolimits_{k = 0}^{ + \infty } {\frac{{{f^{\left( k \right)}}\left( a \right)}}{{\Gamma \left( {k - {\alpha _1} - {\alpha _2}+ 1} \right)}}{{\left( {t - a} \right)}^{k - {\alpha _1} - {\alpha _2}}}} \\
 =&\hspace{-6pt}  {}_a^{\rm{R}}{\mathscr D}_t^{{\alpha _1} + {\alpha _2}}f\left( t \right)\\
 =&\hspace{-6pt}  {}_a^{\rm{R}}{\mathscr D}_t^{{\alpha _2}}{}_a^{\rm{R}}{\mathscr D}_t^{{\alpha _1}}f\left( t \right).
\end{array}
\end{eqnarray}

Let ${\alpha _1},~{\alpha _2}\in \left( {n-1,n} \right)$ and $n\in\mathbb{N}_+$, one has
\begin{eqnarray}\label{Eq19}
\begin{array}{rl}
{}_a^{\rm{C}}{\mathscr D}_t^{{\alpha _1}}{}_a^{\rm{C}}{\mathscr D}_t^{{\alpha _2}}f\left( t \right) =&\hspace{-6pt}  {}_a^{\rm{C}}{\mathscr D}_t^{{\alpha _1}}\sum\nolimits_{k = n}^{ + \infty } {\frac{{{f^{\left( k \right)}}\left( a \right)}}{{\Gamma \left( {k - {\alpha _2} + 1} \right)}}{{\left( {t - a} \right)}^{k - {\alpha _2}}}} \\
 =&\hspace{-6pt}  {}_a^{\rm R}{\mathscr I}_t^{n-{\alpha _1}}\sum\nolimits_{k = n}^{ + \infty } {\frac{{{f^{\left( k \right)}}\left( a \right)}}{{\Gamma \left( {k-n - {\alpha _2}+1} \right)}}{{\left( {t - a} \right)}^{k - {\alpha _2}-n}}} \\
 = &\hspace{-6pt}\sum\nolimits_{k = n}^{ + \infty } {\frac{{{f^{\left( k \right)}}\left( a \right)}}{{\Gamma \left( {k - {\alpha _1} - {\alpha _2} + 1} \right)}}{{\left( {t - a} \right)}^{k - {\alpha _1} - {\alpha _2}}}} \\
 = &\hspace{-6pt} {}_a^{\rm{C}}{\mathscr D}_t^{{\alpha _2}}{}_a^{\rm{C}}{\mathscr D}_t^{{\alpha _1}}f\left( t \right).
\end{array}
\end{eqnarray}
Furthermore, if $n=1$ and ${\alpha _1} + {\alpha _2} \in \left( {0,1} \right)$, one has
\begin{eqnarray}\label{Eq20}
{}_a^{\rm{C}}{\mathscr D}_t^{{\alpha _1}}{}_a^{\rm{C}}{\mathscr D}_t^{{\alpha _2}}f\left( t \right) ={}_a^{\rm{C}}{\mathscr D}_t^{{\alpha _2}}{}_a^{\rm{C}}{\mathscr D}_t^{{\alpha _1}}f\left( t \right) =  {}_a^{\rm{C}}{\mathscr D}_t^{{\alpha _1} + {\alpha _2}}f\left( t \right).
\end{eqnarray}
When ${\beta _k} \in \left( {0,1} \right)$ and $\sum\nolimits_{k = 1}^m {{\beta _k}}  = {\alpha _1} + {\alpha _2}$, one has that
\begin{equation}\label{Eq21}
{}_a^{\rm{C}}{\mathscr D}_t^{{\alpha _1} + {\alpha _2}}f\left( t \right) = {}_a^{\rm{C}}{\mathscr D}_t^{{\beta _1}}{}_a^{\rm{C}}{\mathscr D}_t^{{\beta _2}} \cdots {}_a^{\rm{C}}{\mathscr D}_t^{{\beta _m}}f\left( t \right).
\end{equation}

For any $\alpha\in(n-1,n),~n\le m,~n\in \mathbb{N}_+$ and $m\in\mathbb{N}_+$, one has
\begin{equation}\label{Eq22}
\begin{array}{rl}
\frac{{\rm d}^m}{{\rm d}t^m}f(t) =&\hspace{-6pt}{}_a^{\rm{C}}{\mathscr D}_t^{{\alpha} }{}_a^{\rm{C}}{\mathscr D}_t^{{m-\alpha} }f\left( t \right)\\
=&\hspace{-6pt}{}_a^{\rm{C}}{\mathscr D}_t^{{m-\alpha} }{}_a^{\rm{C}}{\mathscr D}_t^{{\alpha} }f\left( t \right)\\
=&\hspace{-6pt}{}_a^{\rm{R}}{\mathscr D}_t^{{\alpha} }{}_a^{\rm{C}}{\mathscr D}_t^{{m-\alpha} }f\left( t \right)\\
=&\hspace{-6pt}{}_a^{\rm{C}}{\mathscr D}_t^{{m-\alpha} }{}_a^{\rm{R}}{\mathscr D}_t^{{\alpha}}f\left( t \right)\\
=&\hspace{-6pt}{}_a^{\rm{R}}{\mathscr D}_t^{{m-\alpha} }{}_a^{\rm{C}}{\mathscr D}_t^{ {\alpha}}f\left( t \right)\\
=&\hspace{-6pt}{}_a^{\rm{C}}{\mathscr D}_t^{{\alpha} }{}_a^{\rm{R}}{\mathscr D}_t^{{m-\alpha} }f\left( t \right)\\
=&\hspace{-6pt}{}_a^{\rm{R}}{\mathscr D}_t^{{m-\alpha} }{}_a^{\rm{R}}{\mathscr D}_t^{ {\alpha}}f\left( t \right)\\
=&\hspace{-6pt}{}_a^{\rm{R}}{\mathscr D}_t^{{\alpha} }{}_a^{\rm{R}}{\mathscr D}_t^{{m-\alpha} }f\left( t \right).
\end{array}
\end{equation}

With the special assumption on $f(\cdot)$, a handy result in (\ref{Eq18}) is obtained to replace the restricted result in subsection 2.3.6 of \cite{Podlubny:1999Book}. Compared with Theorems 3.3 and 3.4 of \cite{Li:2007AMC}, the proposed results (\ref{Eq19})-(\ref{Eq21}) are more general. Inspired by Property 6 in \cite{Bai:2017ISA}, formula (\ref{Eq22}) provides more effective alternatives.

\begin{itemize}
  \item Case 2: with fractional integral.
\end{itemize}

For any $\alpha\in (n-1,n),~n\in\mathbb{N}_+$ and $\beta>0$, one has
\begin{equation}\label{Eq23}
\begin{array}{rl}
{}_a^{\rm{R}}{\mathscr D}_t^\alpha {}_a^{\rm R}{\mathscr I}_t^\beta f\left( t \right) = &\hspace{-6pt}{}_a^{\rm{R}}{\mathscr D}_t^\alpha \sum\nolimits_{k = 0}^{ + \infty } {\frac{{{f^{\left( k \right)}}\left( a \right)}}{{\Gamma \left( {k + \beta + 1 } \right)}}{{\left( {t - a} \right)}^{k + \beta }}} \\
 =&\hspace{-6pt} \sum\nolimits_{k = 0}^{ + \infty } {\frac{{{f^{\left( k \right)}}\left( a \right)}}{{\Gamma \left( {k + \beta  - \alpha+ 1 } \right)}}{{\left( {t - a} \right)}^{k + \beta  - \alpha }}}\\
 =&\hspace{-6pt} {}_a^{\rm R}{\mathscr I}_t^\beta{}_a^{\rm{R}}{\mathscr D}_t^\alpha f\left( t \right)\\
 =&\hspace{-6pt} {}_a^{\rm{R}}{\mathscr D}_t^{\alpha-\beta} f\left( t \right),
\end{array}
\end{equation}
which implies
\begin{equation}\label{Eq24}
{\textstyle{}_a^{\rm R}{\mathscr D}_t^\alpha {}_a^{\rm{R}}{\mathscr I}_t^\alpha f\left( t \right) ={}_a^{\rm R}{\mathscr I}_t^\alpha {}_a^{\rm{R}}{\mathscr D}_t^\alpha f\left( t \right) =  f\left( t \right).}
\end{equation}

Similarly, formulas (\ref{Eq4}) and (\ref{Eq2}) give
\begin{equation}\label{Eq25}
\begin{array}{rl}
{}_a^{\rm{C}}{\mathscr D}_t^\alpha {}_a^{\rm R}{\mathscr I}_t^\beta f\left( t \right) = &\hspace{-6pt}{}_a^{\rm{C}}{\mathscr D}_t^\alpha \sum\nolimits_{k = 0}^{ + \infty } {\frac{{{f^{\left( k \right)}}\left( a \right)}}{{\Gamma \left( {k + \beta + 1 } \right)}}{{\left( {t - a} \right)}^{k + \beta }}} \\
 =&\hspace{-6pt} \sum\nolimits_{k = 0}^{ + \infty } {\frac{{{f^{\left( k \right)}}\left( a \right)}}{{\Gamma \left( {k + \beta  - \alpha+ 1 } \right)}}{{\left( {t - a} \right)}^{k + \beta  - \alpha }}};
\end{array}
\end{equation}
\begin{equation}\label{Eq26}
\begin{array}{rl}
{}_a^{\rm R}{\mathscr I}_t^\beta {}_a^{\rm{C}}{\mathscr D}_t^\alpha f\left( t \right) =&\hspace{-6pt} {}_a^{\rm R}{\mathscr I}_t^\beta \sum\nolimits_{k = n}^{ + \infty } {\frac{{{f^{\left( k \right)}}\left( a \right)}}{{\Gamma \left( {k - \alpha+ 1 } \right)}}{{\left( {t - a} \right)}^{k - \alpha }}} \\
 =&\hspace{-6pt} \sum\nolimits_{k = n}^{ + \infty } {\frac{{{f^{\left( k \right)}}\left( a \right)}}{{\Gamma \left( {k + \beta  - \alpha+ 1 } \right)}}{{\left( {t - a} \right)}^{k + \beta  - \alpha }}},
\end{array}
\end{equation}
which imply the following special cases
\begin{equation}\label{Eq27}
{\textstyle{}_a^{\rm{C}}{\mathscr D}_t^\alpha {}_a^{\rm R}{\mathscr I}_t^\alpha f\left( t \right) = f\left( t \right);}
\end{equation}
\begin{equation}\label{Eq28}
{\textstyle{}_a^{\rm R}{\mathscr I}_t^\alpha {}_a^{\rm{C}}{\mathscr D}_t^\alpha f\left( t \right) = f\left( t \right) - \sum\nolimits_{k = 0}^{n - 1} {\frac{{{f^{\left( k \right)}}\left( a \right)}}{{\Gamma \left( {k + 1} \right)}}{{\left( {t - a} \right)}^k}}.}
\end{equation}
From definitions, the left hand of (\ref{Eq26}) can be expressed as
\begin{equation}\label{Eq29}
{\textstyle{}_a^{\rm R}{\mathscr I}_t^\alpha {}_a^{\rm{C}}{\mathscr D}_t^\alpha f\left( t \right) = \frac{1}{{\Gamma \left( \alpha  \right)}}\int_a^t {{{\left( {t - \tau } \right)}^{\alpha  - 1}}{}_a^{\rm{C}}{\mathscr D}_t^\alpha } {f}\left( \tau  \right){\rm{d}}\tau .}
\end{equation}
By applying the integral mean value theorem, there exists a variable $\xi\in(a,t)$ such that
\begin{equation}\label{Eq30}
{\textstyle f\left( t \right) = \sum\nolimits_{k = 0}^{n - 1} {\frac{{{f^{\left( k \right)}}\left( a \right)}}{{\Gamma \left( {k + 1} \right)}}} {\left( {t - a} \right)^k} + \frac{{{}_a^{\rm{C}}{\mathscr D}_t^\alpha f\left( \xi  \right)}}{{\Gamma \left( {\alpha  + 1} \right)}}{\left( {t - a} \right)^\alpha }.}
\end{equation}
Note that the newly derived formula (\ref{Eq30}) is similar to the existing results (see (1.5) in \cite{Trujillo:1999JMAA}, Corollary 7 in \cite{Munkhammar:2004Report} and (3.12) in \cite{Odibat:2007AMC}), while they are essentially different.

\subsection{The differentiability on the order}
Assume $f(\cdot)$ can be expressed as a Taylor series at the initial instant, i.e., $f(t)=\sum\nolimits_{k = 0}^{ + \infty } {\frac{{{f^{\left( k \right)}}\left( a \right)}}{{k!}}{{\left( {t - a} \right)}^k}} $, $t\in[a,b]$. Inspired by the discussion when the order tends to the edge of the interval \cite{Podlubny:1999Book,Li:2007AMC}, the continuity of fractional derivatives with respect to the order $\alpha\in(n-1,n),~n\in\mathbb{N}_+$ can be analyzed. For any constant $\nu\in(n-1,n)$, it follows from the continuity of the Gamma function that

\begin{equation}\label{Eq34}
{\textstyle\begin{array}{rl}
\mathop {\lim }\limits_{\alpha  \to \nu - } {}_a^{\rm{R}}{\mathscr D}_t^\alpha f\left( t \right)=&\hspace{-6pt}\mathop {\lim }\limits_{\alpha  \to \nu + } {}_a^{\rm{R}}{\mathscr D}_t^\alpha f\left( t \right)\\
=&\hspace{-6pt} \sum\nolimits_{k = 0}^{{\rm{ + }}\infty } {\frac{{{f^{\left( k \right)}}\left( a \right)}}{{\Gamma \left( {k - \nu+ 1} \right)}}{{\left( {t - a} \right)}^{k - \nu}}}\\
=&\hspace{-6pt} {}_a^{\rm{R}}{\mathscr D}_t^\nu f\left( t \right);
\end{array}}
\end{equation}
\begin{equation}\label{Eq35}
{\textstyle\begin{array}{rl}
\mathop {\lim }\limits_{\alpha  \to \nu - } {}_a^{\rm{C}}{\mathscr D}_t^\alpha f\left( t \right)=&\hspace{-6pt}\mathop {\lim }\limits_{\alpha  \to \nu + } {}_a^{\rm{C}}{\mathscr D}_t^\alpha f\left( t \right)\\
=&\hspace{-6pt} \sum\nolimits_{k = n}^{{\rm{ + }}\infty } {\frac{{{f^{\left( k \right)}}\left( a \right)}}{{\Gamma \left( {k - \nu+ 1} \right)}}{{\left( {t - a} \right)}^{k - \nu}}}\\
=&\hspace{-6pt} {}_a^{\rm{C}}{\mathscr D}_t^\nu f\left( t \right),
\end{array}}
\end{equation}
which means that both the Riemann--Liouville fractional derivative and the Caputo fractional derivative are continuous with respect to their derivative order on the interval $(n-1,n),~n\in\mathbb{N}_+$.

Now, the continuity near the boundaries can be obtained successfully
\begin{equation}\label{Eq36}
{\textstyle\begin{array}{rl}
\mathop {\lim }\limits_{\alpha  \to n - } {}_a^{\rm{R}}{\mathscr D}_t^\alpha f\left( t \right)
=&\hspace{-6pt} \sum\nolimits_{k = n}^{{\rm{ + }}\infty } {\frac{{{f^{\left( k \right)}}\left( a \right)}}{{\Gamma \left( {k - n+ 1} \right)}}{{\left( {t - a} \right)}^{k - n}}}= \frac{{{{\rm{d}}^n}}}{{{\rm{d}}{t^n}}}f\left( t \right),
\end{array}}
\end{equation}
\begin{equation}\label{Eq37}
\begin{array}{rl}
\mathop {\lim }\limits_{\alpha  \to \left( {n - 1} \right) + } {}_a^{\rm{R}}{\mathscr D}_t^\alpha f\left( t \right)=&\hspace{-6pt} \sum\nolimits_{k = n - 1}^{+\infty } {\frac{{{f^{\left( k \right)}}\left( a \right)}}{{\Gamma \left( {k - n+ 2} \right)}}{{\left( {t - a} \right)}^{k - n + 1}}}=\frac{{{{\rm{d}}^{n - 1}}}}{{{\rm{d}}{t^{n - 1}}}}f\left( t \right);
\end{array}
\end{equation}
\begin{equation}\label{Eq38}
{\textstyle\begin{array}{rl}
\mathop {\lim }\limits_{\alpha  \to n - } {}_a^{\rm{C}}{\mathscr D}_t^\alpha f\left( t \right)=&\hspace{-6pt} \sum\nolimits_{k = n}^{{\rm{ + }}\infty } {\frac{{{f^{\left( k \right)}}\left( a \right)}}{{\Gamma \left( {k - n + 1} \right)}}{{\left( {t - a} \right)}^{k - n}}}=\frac{{{{\rm{d}}^n}}}{{{\rm{d}}{t^n}}}f\left( t \right),
\end{array}}
\end{equation}
\begin{equation}\label{Eq39}
\begin{array}{rl}
\mathop {\lim }\limits_{\alpha  \to \left( {n - 1} \right) + } {}_a^{\rm C}{\mathscr D}_t^\alpha f\left( t \right)=&\hspace{-6pt} \sum\nolimits_{k = n}^{{\rm{ + }}\infty } {\frac{{{f^{\left( k \right)}}\left( a \right)}}{{\Gamma \left( {k - n+ 2} \right)}}{{\left( {t - a} \right)}^{k - n + 1}}}=\frac{{{{\rm{d}}^{n - 1}}}}{{{\rm{d}}{t^{n - 1}}}}f\left( t \right) - {f^{\left( {n - 1} \right)}}\left( a \right),
\end{array}
\end{equation}
where $\Gamma (-k) =\pm \infty$, $\forall k \in \mathbb{N}$ is adopted here.

Equations (\ref{Eq36})-(\ref{Eq39}) clearly indicate that the Riemann--Liouville fractional derivative is continuous with respect to the order $\alpha$ on the interval $[n-1,n],~\forall n\in\mathbb{N}_+$. Obviously, this conclusion still holds in the case that $\alpha\in[0,+\infty)$; The Caputo fractional derivative is continuous with respect to $\alpha$ on the interval $(n-1,n]$ and is generally right discontinuous with respect to $\alpha$ at each integer point $\alpha=n$.

Let $\alpha_1\in (0,+\infty)$ and $\alpha_2\in (n-1,n),~n\in\mathbb{N}_+$. We have
\begin{equation}\label{Eq40}
\begin{array}{rl}
\frac{\partial }{{\partial \alpha_1 }}{}_a^{\rm{R}}{\mathscr D}_t^{\alpha_1} f\left( t \right)
=&\hspace{-6pt}\sum\nolimits_{k = 0}^{ + \infty } {{f^{\left( k \right)}}\left( a \right)\frac{\partial }{{\partial \alpha_1 }}\frac{{{{\left( {t - a} \right)}^{k - \alpha_1 }}}}{{\Gamma \left( {k + 1 - \alpha_1 } \right)}}} \\
=&\hspace{-6pt}\sum\nolimits_{k = 0}^{ + \infty } {\frac{{{f^{\left( k \right)}}\left( a \right)}}{{\Gamma \left( {k - \alpha_1+ 1 } \right)}}{{\left( {t - a} \right)}^{k - \alpha_1 }}\phi(k+1,\alpha_1,t,a)};
\end{array}
\end{equation}
\begin{equation}\label{Eq41}
\begin{array}{l}
\frac{\partial }{{\partial \alpha_2 }}{}_a^{\rm{C}}{\mathscr D}_t^{\alpha_2} f\left( t \right)
= \sum\nolimits_{k = n}^{ + \infty } {\frac{{{f^{\left( k \right)}}\left( a \right)}}{{\Gamma \left( {k - \alpha_2+ 1 } \right)}}{{\left( {t - a} \right)}^{k - \alpha_2 }}} \phi(k+1,\alpha_2,t,a),
\end{array}
\end{equation}
where $\phi(k,\alpha,t,a)= {\frac{{\Gamma '\left( {k - \alpha } \right)}}{{\Gamma \left( {k - \alpha } \right)}}  - \ln \left( {t - a} \right)}$ and ${\Gamma '}\left( x \right) = \int_0^{ + \infty } {{{\rm{e}}^{ - t}}{t^{x - 1}}{\rm{ln}}}\left( t \right){\rm{d}}t$. Such a formula can be used for identifying the order.

Alternatively, from Definition \ref{Definition 2} and Definition \ref{Definition 3}, one has
\begin{equation}\label{Eq42}
\hspace{-2pt}\begin{array}{l}
\frac{\partial }{{\partial \alpha_1 }}{}_a^{\rm{R}}{\mathscr D}_t^{\alpha_1} f\left( t \right) =\frac{1}{{\Gamma \left( {n - \alpha_1 } \right)}}\frac{{{{\rm{d}}^n}}}{{{\rm{d}}{t^n}}}\int_a^t {{{\left( {t - \tau } \right)}^{n - \alpha_1  - 1}}f\left( \tau  \right)\phi(n,\alpha_1,t,\tau){\rm{d}}\tau };
\end{array}
\end{equation}
\begin{equation}\label{Eq43}
\begin{array}{l}
\frac{\partial }{{\partial \alpha_2 }}{}_a^{\rm{C}}{\mathscr D}_t^{\alpha_2} f\left( t \right) = \frac{1}{{\Gamma \left( {n - \alpha_2 } \right)}}\int_a^t {{{\left( {t - \tau } \right)}^{n - \alpha_2  - 1}}{f^{\left( n \right)}}\left( \tau  \right)\phi(n,\alpha_2,t,\tau){\rm{d}}\tau }.
\end{array}
\end{equation}

Considering that the initial conditions of the Caputo fractional derivative are irrelevant to the order $\alpha$, one has
\begin{equation}\label{Eq44}
\begin{array}{rl}
\frac{\partial }{{\partial \alpha_2 }}{}_a^{\rm{C}}{\mathscr D}_t^{\alpha_2} f\left( t \right)= &\hspace{-6pt}\frac{\partial }{{\partial \alpha_2 }}{{\mathscr L}_a^{ - 1}}\big\{ {{s^{\alpha_2} }F\left( s \right) - \sum\nolimits_{k = 0}^{n - 1} {{s^{\alpha_2  - k - 1}}{f^{\left( k \right)}}\left( a \right)} } \big\}\\
=&\hspace{-6pt}{{\mathscr L}_a^{ - 1}}\big\{ {\frac{\partial }{{\partial \alpha_2 }}\{{s^{\alpha_2} }F\left( s \right)\} - \sum\nolimits_{k = 0}^{n - 1} {\frac{\rm d }{{\rm d {\alpha_2} }}{s^{\alpha_2  - k - 1}}{f^{\left( k \right)}}\left( a \right)} } \big\}\\
=&\hspace{-6pt}{{\mathscr L}_a^{ - 1}}\big\{ {s^{\alpha_2} }\ln \left( s \right)F\left( s \right)+{s^{\alpha_2} }\frac{\partial }{{\partial {\alpha_2} }}F\left( s \right)\big\}\\
&\hspace{-6pt} - {{\mathscr L}_a^{ - 1}}\big\{\sum\nolimits_{k = 0}^{n - 1} {{s^{{\alpha_2}  - k - 1}}\ln \left( s \right){f^{\left( k \right)}}\left( a \right)} \big\}\\
=&\hspace{-6pt}{}_a^{\rm{C}}{\mathscr D}_t^{\alpha_2} f\left( t \right)*{{\mathscr L}_a^{ - 1}}\left\{\ln \left( s \right) \right\}+{{\mathscr L}_a^{ - 1}}\big\{{s^{\alpha_2} }\frac{\partial }{{\partial {\alpha_2} }}F\left( s \right) \big\},
\end{array}
\end{equation}
where $F\left( s \right) \triangleq\int_a^{ + \infty } {{{\rm{e}}^{ - s\left( {t - a} \right)}}f\left( t \right){\rm{d}}t}$ is the Laplace transform of ${f\left( t \right)}$ with the initial instant $a$, the corresponding inverse Laplace transform ${\mathscr L}_a^{ - 1}\{F\left( s \right)\}\triangleq\frac{1}{{2\pi {\rm{j}}}}\int_{\beta  - {\rm{j}}\infty }^{\beta  + {\rm{j}}\infty } {{F\left( s \right)} {{\rm{e}}^{s\left( {t - a} \right)}}{\rm{d}}s} $, and $*$ stands for the convolution operation with the initial instant $a$, i.e., $p\left( t \right) \ast q\left( t \right) \triangleq \int_a^t {p\left( \tau  \right)q\left( {t + a - \tau } \right){\rm{d}}\tau }$.

After discussing the case of fractional derivative, let us focus on the Riemann-Liouville fractional integral case in the sequel. If $\nu$ is a constant satisfying $n-1<\nu<n$, $n\in\mathbb{N}_+$, one has
\begin{equation}\label{Eq45}
{\textstyle\begin{array}{l}
\mathop {\lim }\limits_{\alpha  \to \nu} {}_a^{\rm{R}}{\mathscr I}_t^\alpha f\left( t \right)= \sum\nolimits_{k = 0}^{+\infty } {\frac{{{f^{\left( k \right)}}\left( a \right)}}{{\Gamma \left( {k + \nu + 1} \right)}}{{\left( {t - a} \right)}^{k + \nu}}}
= {}_a^{\rm{R}}{\mathscr I}_t^\nu f\left( t \right),
\end{array}}
\end{equation}
\begin{equation}\label{Eq46}
{\textstyle\begin{array}{l}
\mathop {\lim }\limits_{\alpha  \to n - } {}_a^{\rm R}{\mathscr I}_t^\alpha f\left( t \right)=\sum\nolimits_{k = 0}^{{\rm{ + }}\infty } {\frac{{{f^{\left( k \right)}}\left( a \right)}}{{\Gamma \left( {k + n + 1} \right)}}{{\left( {t - a} \right)}^{k + n}}}={}_a^{\rm{R}}{\mathscr I}_t^n f\left( t \right),
\end{array}}
\end{equation}
\begin{eqnarray}\label{Eq47}
\begin{array}{l}
\mathop {\lim }\limits_{\alpha  \to \left( {n - 1} \right) + } {}_a^{\rm R}{\mathscr I}_t^\alpha f\left( t \right)= \sum\nolimits_{k = 0}^{+\infty } {\frac{{{f^{\left( k \right)}}\left( a \right)}}{{\Gamma \left( {k + n} \right)}}{{\left( {t - a} \right)}^{k + n - 1}}}={}_a^{\rm{R}}{\mathscr I}_t^{n-1}f\left( t \right),
\end{array}
\end{eqnarray}
which mean that the Riemann--Liouville fractional integral is continuous with respect to the order $\alpha$ on the interval $[n-1,n],~\forall n\in\mathbb{N}_+$ (or, to be more exact, $[0,+\infty)$). In fact, the Riemann--Liouville fractional integral is differentiable with respect to the order $\alpha$ and the corresponding differential can be obtained from different angles as
\begin{equation}\label{Eq48}
\begin{array}{rl}
\frac{\partial }{{\partial \alpha }}{}_a^{\rm R}{\mathscr I}_t^\alpha f\left( t \right)
=&\hspace{-6pt}\sum\nolimits_{k = 0}^{ + \infty } {{f^{\left( k \right)}}\left( a \right)\frac{\partial }{{\partial \alpha }}\frac{{{{\left( {t - a} \right)}^{k + \alpha }}}}{{\Gamma \left( {k + \alpha + 1 } \right)}}} \\
=&\hspace{-6pt}-\sum\nolimits_{k = 0}^{ + \infty } {\frac{{{f^{\left( k \right)}}\left( a \right)}}{{\Gamma \left( {k + \alpha+ 1 } \right)}}{{\left( {t - a} \right)}^{k + \alpha }}\phi(k+1,-\alpha,t,a)};
\end{array}
\end{equation}
\begin{equation}\label{Eq49}
\begin{array}{l}
\frac{\partial }{{\partial \alpha }}{}_a^{\rm R}{\mathscr I}_t^\alpha f\left( t \right) =-\frac{1}{{\Gamma \left( {\alpha } \right)}}\int_a^t {{{\left( {t - \tau } \right)}^{ \alpha  - 1}}f\left( \tau  \right)\phi(0,-\alpha,t,\tau){\rm{d}}\tau };
\end{array}
\end{equation}
\begin{equation}\label{Eq50}
\hspace{-1pt}\begin{array}{rl}
\frac{\partial }{{\partial \alpha }}{}_a^{\rm R}{\mathscr I}_t^\alpha f\left( t \right)=&\hspace{-6pt} \frac{\partial }{{\partial \alpha }}{{\mathscr L}_a^{ - 1}}\{ {{s^{-\alpha} }F\left( s \right)  } \}\\
=&\hspace{-6pt}{{\mathscr L}_a^{ - 1}}\{ {\frac{\partial }{{\partial \alpha }}{s^{-\alpha} }F\left( s \right)} \}\\
=&\hspace{-6pt}{{\mathscr L}_a^{ - 1}}\{ -{s^{-\alpha} }\ln \left( s \right)F\left( s \right)+{s^{-\alpha} }\frac{\partial }{{\partial \alpha }}F\left( s \right)\}\\
=&\hspace{-6pt}-{}_a^{\rm R}{\mathscr I}_t^\alpha f\left( t \right)*{{\mathscr L}_a^{ - 1}}\left\{\ln \left( s \right) \right\}+{{\mathscr L}_a^{ - 1}}\left\{{s^{-\alpha} }\frac{\partial }{{\partial \alpha }}F\left( s \right) \right\}.
\end{array}
\end{equation}

Compared with the related work in \cite{Podlubny:1999Book,Li:2007AMC}, a more detailed investigation is made. Several practical results on the continuity and the differentiability with respect to the order are provided.
\subsection{Scale transform}
Assume $f(\cdot)$ can be expressed as a Taylor series at the initial instant, i.e., $f(t)=\sum\nolimits_{k = 0}^{ + \infty } {\frac{{{f^{\left( k \right)}}\left( a \right)}}{{k!}}{{\left( {t - a} \right)}^k}} $, $t\in[a,b]$. By a scale transform of a function $f(t)$ with respect to an initial instant $a$, it means to replace $f(t)$ by $f(x)$ with an initial instant $c$, where $x = \lambda t + \gamma ,~c = \lambda a + \gamma $, $\lambda$ and $\gamma$ are scaling factors. In the framework of the proposed series, an effective procedure will be constructed to evaluate the effect on fractional derivative operation.

Let $\alpha\in(n-1,n),~n\in \mathbb{N}_+$. One has the following two equations
\begin{equation}\label{Eq51}
{\textstyle {}_a^{\rm{R}}{\mathscr D}_t^\alpha f\left( x \right) = \sum\nolimits_{k = 0}^{ + \infty } {\frac{{{\lambda ^k}{f^{\left( k \right)}}\left( c \right)}}{{\Gamma \left( {k - \alpha+ 1 } \right)}}} {\left( {t - a} \right)^{k - \alpha }};}
\end{equation}
\begin{equation}\label{Eq52}
{\textstyle {}_a^{\rm{C}}{\mathscr D}_t^\alpha f\left( x \right) = \sum\nolimits_{k = n}^{ + \infty } {\frac{{{\lambda ^k}{f^{\left( k \right)}}\left( c \right)}}{{\Gamma \left( {k - \alpha+ 1 } \right)}}} {\left( {t - a} \right)^{k - \alpha }},}
\end{equation}
where the fact $\frac{{{{\rm{d}}^k}}}{{{\rm{d}}{t^k}}}f\left( x \right) = {\lambda ^k}{f^{\left( k \right)}}\left( x \right)={\lambda ^k}\frac{{{{\rm{d}}^k}}}{{{\rm{d}}{x^k}}}f\left( x \right)$ is adopted here. If the differential variable $t$ and the initial instant $a$ are replaced by $x$ and $c$ respectively, one has
\begin{equation}\label{Eq53}
{\textstyle {}_c^{\rm{R}}{\mathscr D}_x^\alpha f\left( x \right) = \sum\nolimits_{k = 0}^{ + \infty } {\frac{{{f^{\left( k \right)}}\left( c \right)}}{{\Gamma \left( {k - \alpha+ 1 } \right)}}} {\left( {x - c} \right)^{k - \alpha }};}
\end{equation}
\begin{equation}\label{Eq54}
{\textstyle {}_c^{\rm{C}}{\mathscr D}_x^\alpha f\left( x \right) = \sum\nolimits_{k = n}^{ + \infty } {\frac{{{f^{\left( k \right)}}\left( c \right)}}{{\Gamma \left( {k - \alpha + 1 } \right)}}} {\left( {x - c} \right)^{k - \alpha }}.}
\end{equation}
Combining (\ref{Eq51})-(\ref{Eq54}) yields
\begin{equation}\label{Eq55}
{}_a^{\rm{R}}{\mathscr D}_t^\alpha f\left( x \right) = {\lambda ^\alpha }{}_c^{\rm{R}}{\mathscr D}_x^\alpha f\left( x \right);
\end{equation}
\begin{equation}\label{Eq56}
{}_a^{\rm{C}}{\mathscr D}_t^\alpha f\left( x \right) = {\lambda ^\alpha }{}_c^{\rm{C}}{\mathscr D}_{x}^\alpha f\left( x \right).
\end{equation}

In (\ref{Eq55}) and (\ref{Eq56}), if one sets $a = 0,~\gamma  = 0,~\lambda  \in \mathbb{C}$ and $\lambda\ne 0$, then $x=\lambda t,~c=a$ and the scale transform operation reduces to a simple multiplication by a constant, namely,
\begin{equation}\label{Eq57}
{}_0^{\rm{R}}{\mathscr D}_t^\alpha f\left( \lambda t \right) = {\lambda ^\alpha }{}_0^{\rm{R}}{\mathscr D}_{\lambda t}^\alpha f\left( \lambda t \right);
\end{equation}
\begin{equation}\label{Eq58}
{}_0^{\rm{C}}{\mathscr D}_t^\alpha f\left( \lambda t \right) = {\lambda ^\alpha }{}_0^{\rm{C}}{\mathscr D}_{\lambda t}^\alpha f\left( \lambda t \right).
\end{equation}

By setting $a =  - \infty ,~\gamma  = 0$ and $\lambda  > 0$, (\ref{Eq55}) and (\ref{Eq56}) could also be simplified as follows
\begin{equation}\label{Eq59}
{}_{-\infty}^{\hspace{6pt}\rm{R}}{\mathscr D}_t^\alpha f\left( \lambda t \right) = {\lambda ^\alpha }{}_{-\infty}^{\hspace{6pt}\rm{R}}{\mathscr D}_{\lambda t}^\alpha f\left( \lambda t \right);
\end{equation}
\begin{equation}\label{Eq60}
{}_{-\infty}^{\hspace{6pt}\rm{C}}{\mathscr D}_t^\alpha f\left( \lambda t \right) = {\lambda ^\alpha }{}_{-\infty}^{\hspace{6pt}\rm{C}}{\mathscr D}_{\lambda t}^\alpha f\left( \lambda t \right),
\end{equation}
respectively.

By setting $\lambda = 1,~\gamma  = -c$ and $g(t)=f(t-a) $, the following time shift property appears from (\ref{Eq55}) and (\ref{Eq56}) as follows
\begin{equation}\label{Eq61}
{}_a^{\rm{R}}{\mathscr D}_t^\alpha g\left( t \right) = {}_0^{\rm{R}}{\mathscr D}_{t - a}^\alpha g\left( t \right);
\end{equation}
\begin{equation}\label{Eq62}
{}_a^{\rm{C}}{\mathscr D}_t^\alpha g\left( t \right) = {}_0^{\rm{C}}{\mathscr D}_{t - a}^\alpha g\left( t \right),
\end{equation}
respectively.

\begin{remark}\label{Remark 2}
In subsection 5.4 of \cite{Oldham:1974Book}, they focus on the scale transform with the same initial instant $a$ and different variable $\lambda t -\lambda a+a$. In Lemma 2.3 of \cite{Duan:2016JCP}, the conversions on both $t$ and $a$ are considered, which are equal to formulas (\ref{Eq55})-(\ref{Eq56}). In addition to providing a concise derivation process, the proposed series enriches the property of the scale transform. Several special and practical cases in (\ref{Eq57})-(\ref{Eq62}) could lead to more useful applications.
\end{remark}

\subsection{The effect of the initial instant}
Let us look into the behaviour of fractional calculus near the initial instant, i.e., $t\to a+$. From formulas (\ref{Eq2}), (\ref{Eq4}) and (\ref{Eq5}), we have
\begin{eqnarray}\label{Eq63}
{\textstyle\begin{array}{rl}
\mathop {\lim }\limits_{t \to a + } {}_a^{\rm{R}}{\mathscr D}_t^\alpha f\left( t \right)=&\hspace{-6pt}\mathop {\lim }\limits_{t \to a + } \frac{{f\left( a \right)}}{{\Gamma \left( {1 - \alpha } \right)}}{\left( {t - a} \right)^{ - \alpha }}\\
=&\hspace{-6pt} \left\{ \begin{array}{ll}
0&\hspace{-6pt},~\alpha  < 0,~f\left( a \right)\neq \infty,\\
f\left( a \right)&\hspace{-6pt},~\alpha  = 0,\\
\pm\infty&\hspace{-6pt},~\alpha  > 0,~\alpha  \notin \mathbb{N}~{\rm{and}}~f\left( a \right)\neq0;
\end{array}\right.
\end{array}}
\end{eqnarray}
\begin{eqnarray}\label{Eq64}
{\textstyle \begin{array}{rl}
\mathop {\lim }\limits_{t \to a + } {}_a^{\rm{C}}{\mathscr D}_t^\alpha f\left( t \right) =&\hspace{-6pt} \mathop {\lim }\limits_{t \to a + } \frac{{{f^{\left( n \right)}}\left( a \right)}}{{\Gamma \left( {n - \alpha + 1 } \right)}}{\left( {t - a} \right)^{n - \alpha }} \\
=&\hspace{-6pt}0,\alpha\in(n-1,n),~n\in\mathbb{N}_+~{\rm{and}}~{f^{\left( n \right)}}\left( a \right) \neq \infty.
\end{array}}
\end{eqnarray}
where $f(\cdot)$ is assumed to be an analytic function in $(a,b)$.

Interestingly, the preceding equations show that when $f\left( a \right)$ is nonzero, its Riemann--Liouville fractional derivative is singular at the initial instant with $\alpha  > 0$ and $\alpha  \notin \mathbb{N}$, while the Caputo one equals to $0$ for finite ${f^{\left( n \right)}}\left( a \right)$. This is why the Riemann--Liouville fractional derivative has strong singularity while the Caputo one performs weak singularity. Compared with the existing work in \cite{Oldham:1974Book,Podlubny:1999Book}, three critical prerequisites were newly provided for (\ref{Eq63}), namely, $f\left( a \right)\neq \infty$, $\alpha  \notin \mathbb{N}$ and $f\left( a \right)\neq0$. The property in (\ref{Eq64}) was originally proposed.

We shift our attention to the behaviour of fractional derivative far from the initial instant, i.e., $t-a\to +\infty$. Recalling the relationship in (\ref{Eq9}), when $\alpha\in(n-1,n),~n\in\mathbb{N}_+$ and $t-a\to +\infty$, $\sum\nolimits_{k = 0}^{n - 1} {\frac{{{f^{\left( k \right)}}\left( a \right)}}{{\Gamma \left( {k + 1 - \alpha } \right)}}{{\left( {t - a} \right)}^{k - \alpha }}}  = 0$ is valid, and therefore
\begin{equation}\label{Eq65}
\mathop {\lim }\limits_{t - a \to  + \infty } {}_a^{\rm{R}}{\mathscr D}_t^\alpha f\left( t \right) = \mathop {\lim }\limits_{t - a \to  + \infty } {}_a^{\rm{C}}{\mathscr D}_t^\alpha f\left( t \right).
\end{equation}
This properly allows us to continue our discussion merely under Riemann--Liouville definition in the sequel.

To further study this issue, the expansion at the current time instead of the initial instant is considered. Similarly, (\ref{Eq65}) becomes
\begin{equation}\label{Eq66}
{\textstyle{}_a^{\rm{R}}{\mathscr D}_t^\alpha f\left( t \right) = \sum\nolimits_{k = n}^{ + \infty } {\left( {\begin{smallmatrix}
\alpha \\
k
\end{smallmatrix}} \right)\frac{{{f^{\left( k \right)}}\left( t \right)}}{{\Gamma \left( {k - \alpha+ 1 } \right)}}{{\left( {t - a} \right)}^{k - \alpha }}},}
\end{equation}
as $t-a\to +\infty$. The characteristics of Riemann--Liouville fractional derivative with three different cases are discussed in the subsequence.

\begin{itemize}
  \item Case 1: $a$ is finite, $t$ is sufficiently positive large.
\end{itemize}

In this case, $ t  \gg \left| a \right|$. One has
\begin{equation}\label{Eq67}
\begin{array}{rl}
{\left( {t - a} \right)^{k - \alpha }} =&\hspace{-6pt} {t^{k - \alpha }}{\left( {1 - \frac{a}{t}} \right)^{k - \alpha }}\\
 =&\hspace{-6pt} {t^{k - \alpha }}\big[ {1 - \left( {k - \alpha } \right)\frac{a}{t} + O\big( {\frac{{{a^2}}}{{{t^2}}}} \big)} \big]\vspace{2pt}\\
 \approx&\hspace{-6pt} {t^{k - \alpha }} - a\left( {k - \alpha } \right){t^{k - \alpha  - 1}}.
\end{array}
\end{equation}

The relationship
\begin{equation}\label{Eq68}
{\textstyle\begin{array}{rl}
\big( {\begin{smallmatrix}
\alpha \\
k
\end{smallmatrix}} \big)\frac{1}{{\Gamma \left( {k - \alpha  + 1} \right)}} =&\hspace{-6pt} \frac{{\Gamma \left( {\alpha  + 1} \right)}}{{\Gamma \left( {k + 1} \right)\Gamma \left( {\alpha  - k + 1} \right)\Gamma \left( {k - \alpha  + 1} \right)}}\\
 =&\hspace{-6pt} \frac{{\Gamma \left( {\alpha  + 1} \right)}}{{k!\left( {\alpha  - k} \right)\Gamma \left( {\alpha  - k} \right)\Gamma \left( {k - \alpha  + 1} \right)}}\\
 =&\hspace{-6pt} \frac{{\Gamma \left( {\alpha  + 1} \right)\sin \left( {\alpha \pi  - k\pi } \right)}}{{k!\left( {\alpha  - k} \right)\pi }}\\
 =&\hspace{-6pt} \frac{{\Gamma \left( {\alpha  + 1} \right)}}{{k!\left( {\alpha  - k} \right)}}\frac{{\sin \left( {\alpha \pi } \right)}}{\pi }{\left( { - 1} \right)^k}
\end{array}}
\end{equation}
follows from the reflection formula.

Substituting (\ref{Eq67}) and (\ref{Eq68}) into (\ref{Eq66}), one obtains
\begin{eqnarray}\label{Eq69}
\begin{array}{rl}
{}_a^{\rm{R}}{\mathscr D}_t^\alpha f\left( t \right)\approx&\hspace{-6pt} {}_0^{\rm{R}}{\mathscr D}_t^\alpha f\left( t \right) - \sum\nolimits_{k = 0}^{ + \infty } {\left( {\begin{smallmatrix}
\alpha \\
k
\end{smallmatrix}} \right)\frac{{{f^{\left( k \right)}}\left( t \right)}}{{\Gamma \left( {k - \alpha+ 1 } \right)}}a\left( {k - \alpha } \right){t^{k - \alpha  - 1}}} \\
\approx&\hspace{-6pt}{}_0^{\rm{R}}{\mathscr D}_t^\alpha f\left( t \right) + \frac{{a\Gamma \left( {\alpha  + 1} \right)\sin \left( {\alpha \pi } \right)}}{{\pi {t^{\alpha  + 1}}}}\sum\nolimits_{k = 0}^{ + \infty } {{{\left( { - 1} \right)}^k}\frac{{{f^{\left( k \right)}}\left( t \right)}}{{k!}}{t^k}}\\
=&\hspace{-6pt} {}_0^{\rm{R}}{\mathscr D}_t^\alpha f\left( t \right) + \frac{{a\Gamma \left( {\alpha  + 1} \right)\sin \left( {\alpha \pi } \right)}}{{\pi {t^{\alpha  + 1}}}}f\left( 0 \right)\\
\approx&\hspace{-6pt}{}_0^{\rm{R}}{\mathscr D}_t^\alpha f\left( t \right),
\end{array}
\end{eqnarray}
which indicates that the impact of the initial instant at the dynamic process $f(t)$ vanishes as $t\to +\infty$. Hence, for large $t$ the fractional derivative with the initial instant $t=a$ can be replaced by, for instance, the counterpart with the initial instant $t=0$. It is noted that case 1 is similar to the result given in \cite{Oldham:1974Book,Podlubny:1999Book}, while for the sake of completeness of this work, we make a discussion from a different viewpoint.

\begin{itemize}
  \item Case 2: $a$ is sufficiently negative large,  $t$ is finite.
\end{itemize}

In this case, $ -a  \gg \left| t \right|$ and one has
\begin{equation}\label{Eq70}
\begin{array}{rl}
{\left( {t - a} \right)^{k - \alpha }} =&\hspace{-6pt} {\left( { - a} \right)^{k - \alpha }}{\left( {1 - \frac{t}{a}} \right)^{k - \alpha }}\\
 =&\hspace{-6pt} {\left( { - a} \right)^{k - \alpha }}\big[ {1 - \left( {k - \alpha } \right)\frac{t}{a} + O\big( {\frac{{{t^2}}}{{{a^2}}}} \big)} \big]\\
 \approx&\hspace{-6pt} {\left( { - a} \right)^{k - \alpha }} + t\left( {k - \alpha } \right){\left( { - a} \right)^{k - \alpha  - 1}}.
\end{array}
\end{equation}

Therefore, in a similar way, formula (\ref{Eq70}) becomes
\begin{eqnarray}\label{Eq71}
\begin{array}{rl}
{}_a^{\rm{R}}{\mathscr D}_t^\alpha f\left( t \right) \approx&\hspace{-6pt} {}_{t + a}^{\hspace{6pt}\rm{R}}{\mathscr D}_t^\alpha f\left( t \right)- \sum\nolimits_{k = 0}^{ + \infty } {\left( {\begin{smallmatrix}
\alpha \\
k
\end{smallmatrix}} \right)\frac{{{f^{\left( k \right)}}\left( t \right)}}{{\Gamma \left( {k - \alpha+ 1 } \right)}}t\left( {k - \alpha } \right){\left( { - a} \right)^{k - \alpha  - 1}}} \\
\approx&\hspace{-6pt}{}_{t + a}^{\hspace{6pt}\rm{R}}{\mathscr D}_t^\alpha f\left( t \right) - \frac{{t\Gamma \left( {\alpha  + 1} \right)\sin \left( {\alpha \pi } \right)}}{{\pi {{\left( { - a} \right)}^{\alpha  + 1}}}}\sum\nolimits_{k = 0}^{ + \infty } {\frac{{{f^{\left( k \right)}}\left( t \right)}}{{k!}}{a^k}} \\
 =&\hspace{-6pt} {}_{t + a}^{\hspace{6pt}\rm{R}}{\mathscr D}_t^\alpha f\left( t \right) - \frac{{t\Gamma \left( {\alpha  + 1} \right)\sin \left( {\alpha \pi } \right)}}{{\pi {{\left( { - a} \right)}^{\alpha  + 1}}}}f\left( {t + a} \right)\\
 \approx&\hspace{-6pt}{}_{t + a}^{\hspace{6pt}\rm{R}}{\mathscr D}_t^\alpha f\left( t \right).
\end{array}
\end{eqnarray}
It can be concluded that, under certain conditions on $f(t)$ with large negative value of $a$, the fractional derivative with a fixed initial instant ${}_a^{\rm{R}}{\mathscr D}_t^\alpha f\left( t \right)$ can be replaced by the one with a moving initial instant ${}_{t + a}^{\hspace{6pt}\rm{R}}{\mathscr D}_t^\alpha f\left( t \right)$, whose memory length is fixed rather than being increasing with time. Notably, although a similar result ${}_a^{\rm{R}}{\mathscr D}_t^\alpha f\left( t \right) \approx{}_{t - a}^{\hspace{6pt}\rm{R}}{\mathscr D}_t^\alpha f\left( t \right)$ was proposed in \cite{Podlubny:1999Book}, it is somewhat incomplete. Accordingly, the discussion in case 2 was given.

\begin{itemize}
  \item Case 3: $a$ is sufficiently negative large, $t$ is sufficiently positive large.
\end{itemize}

If $\mathop {\lim }\limits_{\tiny a \to  - \infty , t \to  + \infty} \left| {\frac{a}{t}} \right| < 1$, the situation is similar to case 1 ($\left| {\frac{a}{t}} \right| =  0 $) with
\begin{equation}\label{Eq72}
{}_a^{\rm{R}}{\mathscr D}_t^\alpha f\left( t \right) \approx {}_{0}^{\rm{R}}{\mathscr D}_t^\alpha f\left( t \right).
\end{equation}

If $\mathop {\lim }\limits_{\tiny a \to  - \infty , t \to  + \infty} \left| {\frac{a}{t}} \right| > 1$, the situation is similar to case 2 ($\left| {\frac{a}{t}} \right| =  + \infty $) with
\begin{equation}\label{Eq73}
{}_a^{\rm{R}}{\mathscr D}_t^\alpha f\left( t \right) \approx {}_{t + a}^{\hspace{6pt}\rm{R}}{\mathscr D}_t^\alpha f\left( t \right).
\end{equation}

Take into account the transitional case, i.e., $\mathop {\lim }\limits_{\tiny a \to  - \infty \hfill\atop
\tiny t \to  + \infty \hfill}\left| {\frac{a}{t}} \right| = 1$, which can also be distinctly written as $t + a = 0$, one has
\begin{equation}\label{Eq74}
{}_a^{\rm{R}}{\mathscr D}_t^\alpha f\left( t \right) \approx {}_{t + a}^{\hspace{6pt}\rm{R}}{\mathscr D}_t^\alpha f\left( t \right) = {}_0^{\rm{R}}{\mathscr D}_t^\alpha f\left( t \right).
\end{equation}

Note that case 3 is completely new. In fact, all the conclusions in the above three cases can be reached if the series in (\ref{Eq3}) is adopted instead of (\ref{Eq66}).

\subsection{Derivative computation}
The provided series representation could contribute to fractional derivative calculation in two ways, that is to say, analytical computation and numerical approximation.

\begin{itemize}
  \item Application 1: analytical computation
\end{itemize}

As an example, the Heaviside function $H\left( {t-a } \right)$ occurring at $t=a$ is selected. In view of the convergence domain of infinite series, it is established that
\begin{eqnarray}\label{Eq75}
{\textstyle
\begin{array}{rl}
{}_0^{\rm{R}}{\mathscr D}_t^\alpha H\left( {t-a } \right)=&\hspace{-6pt} \sum\nolimits_{k = 0}^{ + \infty } {\frac{{{H^{\left( k \right)}}\left( { \max\{a,0\}}-a \right)}}{{\Gamma \left( {k - \alpha+ 1 } \right)}}{\left( t-{ \max\{a,0\}}\right)^{k - \alpha }}}\\
=&\hspace{-6pt} \left\{ \begin{array}{rl}
\frac{{{{\left( {t - a} \right)}^{ - \alpha }}}}{{\Gamma \left( {1 - \alpha } \right)}}&\hspace{-6pt},~a > 0,\\
\frac{{{t^{ - \alpha }}}}{{\Gamma \left( {1 - \alpha } \right)}}&\hspace{-6pt},~a \le 0.
\end{array} \right.
\end{array}}
\end{eqnarray}

After applying (\ref{Eq5}), the derivative of the exponential function appears
\begin{equation}\label{Eq76}
\begin{array}{rl}
{}_a^{\rm{R}}{\mathscr D}_t^\alpha {{\rm{e}}^{\lambda t}} =&\hspace{-6pt} \sum\nolimits_{k = 0}^{ + \infty } {\frac{{{\lambda ^k}{{\rm{e}}^{\lambda a}}}}{{\Gamma \left( {k - \alpha+ 1 } \right)}}{{\left( {t - a} \right)}^{k - \alpha }}} \\
=&\hspace{-6pt} {\left( {t - a} \right)^{-\alpha} }{{\mathcal E}_{1,1 - \alpha }}\left( {\lambda t - \lambda a} \right),
\end{array}
\end{equation}
where ${{\mathcal E}_{x,y }}\left( z \right)=\sum\nolimits_{k = 0}^{ + \infty } \frac{z^k}{\Gamma \left( k x + y  \right)}$ is the Mittag--Leffler function. For $n\in\mathbb{N}_+$ and $n-1<\alpha<n$, with the adoption of (\ref{Eq5}), one has
\begin{equation}\label{Eq77}
{}_a^{\rm{C}}{\mathscr D}_t^\alpha {{\rm{e}}^{\lambda t}} ={\left( {t - a} \right)^{n-\alpha} }{{\mathcal E}_{1,1 +n- \alpha }}\left( {\lambda t - \lambda a} \right).
\end{equation}

Moreover, such a derivative can be also expressed as
\begin{equation}\label{Eq78}
\begin{array}{rl}
{}_a^{\rm{R}}{\mathscr D}_t^\alpha {{\rm{e}}^t} =&\hspace{-6pt} \sum\nolimits_{k = 0}^{ + \infty } {\frac{{{{\rm{e}}^a}}}{{\Gamma \left( {k - \alpha+ 1 } \right)}}{{\left( {t - a} \right)}^{k - \alpha }}} \\
 =&\hspace{-6pt} {{\rm{e}}^t}{\left( {t - a} \right)^{ - \alpha }}{{\rm{e}}^{a - t}}\sum\nolimits_{k = 0}^{ + \infty } {\frac{{{{\left( {t - a} \right)}^k}}}{{\Gamma \left( {k - \alpha + 1 } \right)}}} \\
 =&\hspace{-6pt} {{\rm{e}}^t}{g^*}\left( { - \alpha ,t - a} \right),
\end{array}
\end{equation}
where ${g^*}\left( {c,x} \right) = {x^c}{{\rm{e}}^{ - x}}\sum\nolimits_{k = 0}^{ + \infty } {\frac{{{x^k}}}{{\Gamma \left( {k + c+ 1} \right)}}} $ is the auxiliary incomplete Gamma function \cite{Gautschi:1979TMS}. By exploring the asymptotic behavior with a bounded constant $\left| c \right|$, i.e., $\mathop {\lim }\limits_{x \to  + \infty } {g^*}\left( {c,x} \right) = 1$, one has
\begin{equation}\label{Eq79}
\mathop {\lim }\limits_{t - a \to  + \infty } {}_a^{\rm{R}}{\mathscr D}_t^\alpha {{\rm{e}}^t} = {{\rm{e}}^t}.
\end{equation}

Revisiting the scale transform in (\ref{Eq57}), one has that
\begin{equation}\label{Eq80}
{}_0^{\rm{R}}{\mathscr D}_t^\alpha {{\rm{e}}^{\lambda t}} = {\lambda ^\alpha }{{\rm{e}}^{\lambda t}}
\end{equation}
holds for sufficiently large $t$ and $\lambda\in\mathbb{C}$.

It follows from the property in (\ref{Eq59}) that
\begin{equation}\label{Eq81}
{}_{ - \infty }^{\hspace{6pt}\rm{R}}{\mathscr D}_t^\alpha {{\rm{e}}^{\lambda t}} = {\lambda ^\alpha }{{\rm{e}}^{\lambda t}},
\end{equation}
with $\lambda>0$.

A straight calculation on (\ref{Eq80}) yields
\begin{equation}\label{Eq82}
\begin{array}{rl}
{}_{ 0 }^{\rm{R}}{\mathscr D}_t^\alpha {\rm{sin}}\left( {\lambda t + \varphi } \right) =&\hspace{-6pt} \frac{{{{\left( {{\rm{i}}\lambda } \right)}^\alpha }{{\rm{e}}^{{\rm{i}}\left( {\lambda t + \varphi } \right)}} - {{\left( { - {\rm{i}}\lambda } \right)}^\alpha }{{\rm{e}}^{ - {\rm{i}}\left( {\lambda t + \varphi } \right)}}}}{{2{\rm{i}}}}\vspace{2pt}\\
 =&\hspace{-6pt} {\lambda ^\alpha }\frac{{{{\rm{e}}^{{\rm{i}}( {\lambda t + \varphi  + \frac{{\alpha \pi }}{2}} )}} - {{\rm{e}}^{ - {\rm{i}}( {\lambda t + \varphi  + \frac{{\alpha \pi }}{2}} )}}}}{{2{\rm{i}}}}\\
 =&\hspace{-6pt} {\lambda ^\alpha }{\rm{sin}}\left( {\lambda t + \varphi  + \frac{{\alpha \pi }}{2}} \right).
\end{array}
\end{equation}
Likewise, (\ref{Eq81}) leads to
\begin{equation}\label{Eq83}
{\textstyle {}_{ -\infty }^{\hspace{6pt}\rm{R}}{\mathscr D}_t^\alpha {\rm{sin}}\left( {\lambda t + \varphi } \right)= {\lambda ^\alpha }{\rm{sin}}\left( {\lambda t + \varphi  + \frac{{\alpha \pi }}{2}} \right),}
\end{equation}
which coincides with the result on page 311 of \cite{Podlubny:1999Book}. It should point out that the developed series representation is helpful to calculating the fractional derivative of a given function. Compared with the existing results, more specific conditions are provided to deduce the related results.

\begin{itemize}
  \item Application 2: numerical approximation
\end{itemize}

Unlike Application 1, it is generally difficult to access the analytic solution. At this point, it becomes reasonable and practical to consider numerical approximation, that is to say, only limited items will be reserved and the long tail is truncated. Similar to the method in \cite{Atanackovic:2004FCAA}, the following formulas follow from (\ref{Eq2}), (\ref{Eq5}), (\ref{Eq3}) and (\ref{Eq7}) immediately
\begin{equation}\label{Eq85}
{\textstyle {}_a^{\rm{C}}{\mathscr D}_t^\alpha f\left( t \right) \approx \sum\nolimits_{k = n}^N {\frac{{{f^{\left( k \right)}}\left( a \right)}}{{\Gamma \left( {k - \alpha+ 1 } \right)}}{{\left( {t - a} \right)}^{k - \alpha }}},}
\end{equation}
\begin{equation}\label{Eq84}
{\textstyle {}_a^{\rm{R}}{\mathscr D}_t^\alpha f\left( t \right) \approx  \sum\nolimits_{k = 0}^N {\frac{{{f^{\left( k \right)}}\left( a \right)}}{{\Gamma \left( {k - \alpha+ 1 } \right)}}{{\left( {t - a} \right)}^{k - \alpha }}},}
\end{equation}
\begin{equation}\label{Eq87}
{\textstyle{}_a^{\rm{C}}{\mathscr D}_t^\alpha f\left( t \right) \approx\sum\nolimits_{k = n}^N {\left( {\begin{smallmatrix}
{\alpha  - n}\\
{k - n}
\end{smallmatrix}} \right)\frac{{{f^{\left( k \right)}}\left( t \right)}}{{\Gamma \left( {k - \alpha+ 1 } \right)}}{{\left( {t - a} \right)}^{k - \alpha }}},}
\end{equation}
\begin{equation}\label{Eq86}
{\textstyle{}_a^{\rm{R}}{\mathscr D}_t^\alpha f\left( t \right) \approx\sum\nolimits_{k = 0}^N {\left( {\begin{smallmatrix}
\alpha \\
k
\end{smallmatrix}} \right)\frac{{{f^{\left( k \right)}}\left( t \right)}}{{\Gamma \left( {k - \alpha + 1 } \right)}}{{\left( {t - a} \right)}^{k - \alpha }}},}
\end{equation}
where $n-1<\alpha<n$, $n\in\mathbb{N}_+$ and $N\in\mathbb{N}_+$.

To provide an intuitive analysis for the radius of convergence, we use (\ref{Eq2}) as an example. Suppose that the Taylor series expansion of $f(t)$, i.e.,
\begin{equation}\label{Eq88}
{\textstyle  f\left( t \right) =  \sum\nolimits_{k = 0}^{+\infty} {\frac{{{f^{\left( k \right)}}\left( a \right)}}{{\Gamma \left( {k + 1} \right)}}{{\left( {t - a} \right)}^{k  }}}}
\end{equation}
is convergent for $\left| {x - a} \right| < R$ and (\ref{Eq2}) is convergent for $x \in ( {a,a + \hat R} )$. Then,
\begin{equation}\label{Eq89}
{\textstyle \begin{array}{rl}
\hat R =&\hspace{-6pt} \mathop {\lim }\limits_{k \to  + \infty } \big| {\frac{{{f^{\left( k \right)}}\left( a \right)}}{{\Gamma \left( {k- \alpha + 1 } \right)}}\frac{{\Gamma \left( {k - \alpha+ 2 } \right)}}{{{f^{\left( {k + 1} \right)}}\left( a \right)}}} \big|\\
 =&\hspace{-6pt} \mathop {\lim }\limits_{k \to  + \infty } \big| {\frac{{\left( {k - \alpha+ 1 } \right){f^{\left( k \right)}}\left( a \right)}}{{{f^{\left( {k + 1} \right)}}\left( a \right)}}} \big|\\
 =&\hspace{-6pt} \mathop {\lim }\limits_{k \to  + \infty } \big| {\frac{{\left( {k + 1} \right){f^{\left( k \right)}}\left( a \right)}}{{{f^{\left( {k + 1} \right)}}\left( a \right)}}} \big|\\
 =&\hspace{-6pt} \mathop {\lim }\limits_{k \to  + \infty } \big| {\frac{{{f^{\left( k \right)}}\left( a \right)}}{{\Gamma \left( {k + 1} \right)}}\frac{{\Gamma \left( {k + 2} \right)}}{{{f^{\left( {k + 1} \right)}}\left( a \right)}}} \big|\\
 =&\hspace{-6pt} R,
\end{array}}
\end{equation}
which means that the series (\ref{Eq2}) shares the same convergence radius with the series in (\ref{Eq88}).

The convergence of (\ref{Eq88}) implies that, for any given constant $\epsilon  > 0$, there exists an integer $N=N(\epsilon)$, such that
\begin{equation}\label{Eq90}
{\textstyle\big| {\sum\nolimits_{k = N+1}^{ + \infty } {\frac{{{f^{\left( k \right)}}\left( a \right)}}{{\Gamma \left( {k + 1} \right)}}{{\left( {t - a} \right)}^k}} } \big| < \epsilon.}
\end{equation}
Consequently, the truncation error of (\ref{Eq85}) satisfies
\begin{equation}\label{Eq91}
\begin{array}{rl}
\big| {\sum\nolimits_{k = N+1}^{ + \infty } {\frac{{{f^{\left( k \right)}}\left( a \right)}}{{\Gamma \left( {k - \alpha+ 1 } \right)}}{{\left( {t - a} \right)}^{k - \alpha }}} } \big| =&\hspace{-6pt} \big| {{}_a^{\rm{R}}{\mathscr D}_t^\alpha \sum\nolimits_{k = N+1}^{ + \infty } {\frac{{{f^{\left( k \right)}}\left( a \right)}}{{\Gamma \left( {k + 1} \right)}}{{\left( {t - a} \right)}^k}} } \big|\\
 \le&\hspace{-6pt} \big| {{}_a^{\rm R}{\mathscr I}_t^{-\alpha} \big| {\sum\nolimits_{k = N+1}^{ + \infty } {\frac{{{f^{\left( k \right)}}\left( a \right)}}{{\Gamma \left( {k + 1} \right)}}{{\left( {t - a} \right)}^k}} } \big|} \big|\\
 <&\hspace{-6pt} \big| {\frac{\epsilon }{{\Gamma \left( { - \alpha } \right)}}\int_a^t {{{\left( {t - \tau } \right)}^{ - \alpha  - 1}}{\rm{d}}\tau } } \big|\\
 =&\hspace{-6pt} \frac{{\epsilon {{\left( {t - a} \right)}^{ - \alpha }}}}{{\left| {\Gamma \left( {1 - \alpha } \right)} \right|}}.
\end{array}
\end{equation}

Analogously, for any $\epsilon>0$, there exists an integer $N=N(\epsilon)$, such that
\begin{equation}\label{Eq92}
{\textstyle\big| {\sum\nolimits_{k = N+1}^{ + \infty } {\frac{{{f^{\left( k \right)}}\left( t \right)}}{{\Gamma \left( {k + 1} \right)}}{{\left( {\tau - t} \right)}^k}} } \big| < \epsilon.}
\end{equation}

Correspondingly, the truncation error of (\ref{Eq87}) satisfies
\begin{equation}\label{Eq93}
\hspace{-3pt}\begin{array}{rl}
\big|\sum\nolimits_{k = N + 1}^{ + \infty } {\left( {\begin{smallmatrix}
\alpha \\
k
\end{smallmatrix}} \right)\frac{{{f^{\left( k \right)}}\left( t \right)}}{{\Gamma \left( {k - \alpha + 1 } \right)}}{{\left( {t - a} \right)}^{k - \alpha }}} \big| =&\hspace{-6pt} \big|\frac{1}{{\Gamma \left( { - \alpha } \right)}}\int_a^t {\sum\nolimits_{k = N + 1}^{ + \infty } {{{\left( { - 1} \right)}^k}\frac{{{f^{\left( k \right)}}\left( t \right)}}{{\Gamma \left( {k + 1} \right)}}{{\left( {t - \tau } \right)}^{k - \alpha  - 1}}} {\rm{d}}\tau } \big|\\
 \le&\hspace{-6pt}  \big|\frac{1}{{\Gamma \left( { - \alpha } \right)}}\int_a^t {{{\left( {t - \tau } \right)}^{ - \alpha  - 1}}\big| {\sum\nolimits_{k = N + 1}^{ + \infty } {\frac{{{f^{\left( k \right)}}\left( t \right)}}{{\Gamma \left( {k + 1} \right)}}{{\left( {t - \tau } \right)}^k}} } \big|{\rm{d}}\tau } \big|\\
 <&\hspace{-6pt} \big|\frac{\epsilon}{{\Gamma \left( { - \alpha } \right)}}\int_a^t {{{\left( {t - \tau } \right)}^{ - \alpha  - 1}}{\rm{d}}\tau } \big|\\
 =&\hspace{-6pt} \frac{{{\epsilon{\left( {t - a} \right)}^{ - \alpha }}}}{{\left| {\Gamma \left( {1 - \alpha } \right)} \right|}}.
\end{array}\hspace{-6pt}
\end{equation}

Note that the proposed formulas in (\ref{Eq85})-(\ref{Eq86}) can be also be regarded as the approximation of (\ref{Eq2}), (\ref{Eq5}), (\ref{Eq3}) and (\ref{Eq7}) respectively. intuitively, the larger $N$, the smaller the approximation error.

\subsection{Variable order case}
The previous discussions focus on fractional integral or derivative with constant order. To meet the strong demand in practical applications, this subsection investigates the variable order case.

\begin{definition} \label{Definition 4}
The $\alpha\left( t  \right)$-th Riemann--Liouville fractional integral of a function $f(t)$ is defined as
\begin{equation}\label{Eq94}
{\textstyle
{{}_a^{\rm R}}{{\mathscr I}_{t}^{\alpha\left( t  \right)}}{f\left( t \right)}  \triangleq \frac{1}{{\Gamma \left( \alpha\left( t  \right)  \right)}}\int_{{a}}^t {{{\left( {t - \tau } \right)}^{\alpha\left( t  \right) -1}}{ f\left( \tau  \right)}{\rm{d}}\tau },}
\end{equation}
where $\alpha\left( t  \right)>0$ is the integral order and $a$ is the initial instant.
\end{definition}

Considering $\alpha\left( t  \right) > 0$ and $\nu  >  - 1$, then one has
\begin{equation}\label{Eq95}
\begin{array}{rl}
{}_a^{\rm{R}}{\mathscr I}_t^{\alpha\left( t \right)} {\left( {t - a} \right)^\nu } =&\hspace{-6pt} \frac{1}{{\Gamma \left( \alpha\left( t \right)  \right)}}\int_a^t {{{\left( {t - \tau } \right)}^{\alpha\left( t \right)  - 1}}{{\left( {\tau  - a} \right)}^\nu }{\rm{d}}\tau } \\
 =&\hspace{-6pt} \frac{1}{{\Gamma \left( \alpha\left( t \right)  \right)}}\int_0^1 {{{\left( {t - a} \right)}^{\alpha\left( t \right)  - 1}}{{\left( {1 - u} \right)}^{\alpha\left( t \right)  - 1}}{{\left( {t - a} \right)}^\nu }{u^\nu }\left( {t - a} \right){\rm{d}}u} \\
 =&\hspace{-6pt} \frac{{{{\left( {t - a} \right)}^{\alpha\left( t \right)  + \nu }}}}{{\Gamma \left( \alpha\left( t \right)  \right)}}\int_0^1 {{u^\nu }{{\left( {1 - u} \right)}^{\alpha\left( t \right)  - 1}}{\rm{d}}u} \\
 =&\hspace{-6pt} \frac{{{{\left( {t - a} \right)}^{\alpha\left( t \right)  + \nu }}}}{{\Gamma \left( \alpha\left( t \right)  \right)}}B\left( {\nu  + 1,\alpha\left( t \right) } \right)\\
 =&\hspace{-6pt} \frac{{\Gamma \left( {\nu  + 1} \right)}}{{\Gamma \left( {\nu  + \alpha\left( t \right)  + 1} \right)}}{\left( {t - a} \right)^{\nu  + \alpha\left( t \right) }},
\end{array}
\end{equation}
where $\tau  = u\left( {t - a} \right) + a$ and $B\left( {p,q} \right) = \frac{{\Gamma \left( p \right)\Gamma \left( q \right)}}{{\Gamma \left( {p + q} \right)}}$, $p,q>0$ are adopted here.

If $f(t)$ can be expressed as a Taylor series expanded at $t=a$, then its $\alpha\left( t  \right)$-th Riemann--Liouville integral can be rewritten as
\begin{equation}\label{Eq96}
\begin{array}{rl}
{_a^{\rm R}}{\mathscr I}_t^{\alpha\left( t \right)} f\left( t \right) =&\hspace{-6pt} {}_a^{\rm{R}}{\mathscr I}_t^{\alpha\left( t \right)} \sum\nolimits_{k = 0}^{ + \infty } {\frac{{{f^{\left( k \right)}}\left( a \right)}}{{k!}}{{\left( {t - a} \right)}^k}} \\
=&\hspace{-6pt} \sum\nolimits_{k = 0}^{ + \infty } {\frac{{{f^{\left( k \right)}}\left( a \right)}}{{k!}}{}_a^{\rm{R}}{\mathscr I}_t^{\alpha\left( t \right)} {{\left( {t - a} \right)}^k}}\\
=&\hspace{-6pt} \sum\nolimits_{k = 0}^{ + \infty } {\frac{{{f^{\left( k \right)}}\left( a \right)}}{{\Gamma \left( {k + \alpha\left( t \right) + 1} \right)}}{{\left( {t - a} \right)}^{k + \alpha\left( t \right) }}},
\end{array}
\end{equation}
where $\alpha\left( t \right)>0$.

If $f(\tau)$ can be expressed as a Taylor series expanded at $\tau=t$, then another useful representation can be derived as
\begin{equation}\label{Eq97}
\begin{array}{rl}
{}_a^{\rm R}{\mathscr I}_t^{\alpha\left( t \right)} f\left( t \right) =&\hspace{-6pt} \frac{1}{{\Gamma \left( \alpha\left( t \right)  \right)}}\int_a^t {{{\left( {t - \tau } \right)}^{\alpha\left( t \right)  - 1}}f\left( \tau  \right){\rm{d}}\tau } \\
 =&\hspace{-6pt} \frac{1}{{\Gamma \left( \alpha\left( t \right)  \right)}}\int_a^t {\sum\nolimits_{k = 0}^{ + \infty } {{{\left( { - 1} \right)}^k}\frac{{{f^{\left( k \right)}}\left( t \right)}}{{\Gamma \left( {k + 1} \right)}}} {{\left( {t - \tau } \right)}^{k + \alpha\left( t \right)  - 1}}{\rm{d}}\tau } \\
 =&\hspace{-6pt} \sum\nolimits_{k = 0}^{ + \infty } {{{\left( { - 1} \right)}^k}\frac{{{f^{\left( k \right)}}\left( t \right)}}{{\Gamma \left( \alpha\left( t \right)  \right)\Gamma \left( {k + 1} \right)}}} \int_a^t {{{\left( {t - \tau } \right)}^{k + \alpha\left( t \right)  - 1}}{\rm{d}}\tau } \\
 =&\hspace{-6pt} \sum\nolimits_{k = 0}^{ + \infty } {{{\left( { - 1} \right)}^k}\frac{{{f^{\left( k \right)}}\left( t \right)}}{{\Gamma \left( \alpha\left( t \right)  \right)\Gamma \left( {k + 1} \right)\left( {k + \alpha\left( t \right) } \right)}}{{\left( {t - a } \right)}^{k + \alpha\left( t \right) }}} \\
 =&\hspace{-6pt} \sum\nolimits_{k = 0}^{ + \infty } {\big( {\begin{smallmatrix}
{ - \alpha\left( t \right) }\\
k
\end{smallmatrix}} \big)\frac{{{f^{\left( k \right)}}\left( t \right)}}{{\Gamma \left( {k+ \alpha\left( t \right)+ 1} \right)}}{{\left( {t - a} \right)}^{k + \alpha\left( t \right) }}}.
\end{array}
\end{equation}

Similarly, with the help of such a fractional integral in (\ref{Eq94}), the corresponding fractional derivatives could be established, respectively. Let $n\left( t  \right)\in\mathbb{N}_+$ and $n\left( t  \right)-1<\alpha\left( t  \right)<n\left( t  \right)$ for any $t\ge a$.

\begin{definition} \label{Definition 5}
The $\alpha\left( t  \right)$-th Riemann--Liouville fractional derivative of a function $f(t)$ is defined as ${}^{\rm R}_a{\mathscr D}_{t}^{\alpha\left( t  \right)} f\left( t \right) \triangleq \frac{{\rm{d}}^{{n\left( t  \right)}}}{{{\rm{d}}t^{{n\left( t  \right)}}}} {{}_a^{\rm R}{\mathscr I}_{t}^{{n\left( t  \right)}-\alpha\left( t  \right)}} f\left( t \right)$.
\end{definition}
\begin{definition} \label{Definition 6}
The $\alpha\left( t  \right)$-th Caputo fractional derivative of a function $f(t)$ is defined as ${}^{\rm C}_a{\mathscr D}_{t}^{\alpha\left( t  \right)} f\left( t \right) \triangleq  {{}_a^{\rm R}}{\mathscr I}_{t}^{n\left( t  \right)-\alpha\left( t  \right)}\frac{{\rm{d}}^{{n\left( t  \right)}}}{{{\rm{d}}t^{{n\left( t  \right)}}}} f\left( t \right)$.
\end{definition}

Notably, though the introduction of the time varying order extends the application of fractional calculus, many difficulties appear simultaneously. The Laplace transform of variable order fractional integrals/derivatives in Definitions \ref{Definition 4}-\ref{Definition 6} is not concise. Moreover, it is difficult or even impossible to derive the equivalent form in series.

When $n\left( t  \right)$, ${n\left( t  \right)}-\alpha\left( t  \right)$ are not differentiable or even not continuous, extra assumptions should be provided. Assume that $\alpha\left( t  \right)$ is finite and piecewise constant, more specially, $\alpha\left( t  \right)=\alpha_i$, and $n\left( t  \right)=n_i$ hold $t\in[t_{i-1},t_i)$, where $t_0=a$, $t=t_m$. Then, for any $t\in[t_{i-1},t_i)$, one has
\begin{equation}\label{Eq98}
{\textstyle
\begin{array}{rl}
{}^{\rm R}_a{\mathscr D}_{t}^{\alpha\left( t \right)} f\left( t \right) =&\hspace{-6pt}\frac{{\rm{d}}^{n\left( t \right)}}{{{\rm{d}}t^{n\left( t \right)}}} \sum\nolimits_{k = 0}^{ + \infty } {\frac{{{f^{\left( k \right)}}\left( a \right)}}{{\Gamma \left( {k + n\left( t \right)-\alpha\left( t \right)+ 1 } \right)}}{{\left( {t - a} \right)}^{k + n\left( t \right)-\alpha\left( t \right) }}}\\
=&\hspace{-6pt}\frac{{\rm{d}}^{n_i}}{{{\rm{d}}t^{n_i}}} \sum\nolimits_{k = 0}^{ + \infty } {\frac{{{f^{\left( k \right)}}\left( a \right)}}{{\Gamma \left( {k + n_i-\alpha_i+ 1 } \right)}}{{\left( {t - a} \right)}^{k + n_i-\alpha_i }}}\\
=&\hspace{-6pt} \sum\nolimits_{k = 0}^{+\infty}  {\frac{{{f^{\left( k \right)}}\left( a \right)}}{{\Gamma \left( { k- \alpha_i + 1 } \right)}}} {\left( {t - a} \right)^{k - \alpha_i }}.
\end{array}}
\end{equation}
Considering the whole interval $[a,t]$, it follows
\begin{equation}\label{Eq99}
{\textstyle
\begin{array}{l}
{}^{\rm R}_a{\mathscr D}_{t}^{\alpha\left( t \right)} f\left( t \right) =\sum\nolimits_{k = 0}^{+\infty}  {\frac{{{f^{\left( k \right)}}\left( a \right)}}{{\Gamma \left( { k- \alpha\left( t \right) + 1 } \right)}}} {\left( {t - a} \right)^{k - \alpha\left( t \right) }}.
\end{array}}
\end{equation}

Likewise, it follows
\begin{eqnarray}\label{Eq100}
\begin{array}{rl}
{}_a^{\rm{R}}{\mathscr D}_t^{\alpha\left( t \right)} f\left( t \right)=&\hspace{-6pt}\frac{{{{\rm{d}}^{n\left( t \right)}}}}{{{\rm{d}}{t^{n\left( t \right)}}}} \sum\nolimits_{k = 0}^{ + \infty } {\big( {\begin{smallmatrix}
{ \alpha\left( t \right)-n\left( t \right) }\\
k
\end{smallmatrix}} \big)\frac{{{f^{\left( k \right)}}\left( t \right)}}{{\Gamma \left( {k + n\left( t \right)-\alpha\left( t \right)+ 1} \right)}}{{\left( {t - a} \right)}^{k + n\left( t \right)- \alpha\left( t \right) }}} \\
=&\hspace{-6pt}\frac{{{{\rm{d}}^{n_i}}}}{{{\rm{d}}{t^{n_i}}}} \sum\nolimits_{k = 0}^{ + \infty } {\big( {\begin{smallmatrix}
{ \alpha_i-n_i }\\
k
\end{smallmatrix}} \big)\frac{{{f^{\left( k \right)}}\left( t \right)}}{{\Gamma \left( {k + n_i-\alpha_i+ 1} \right)}}{{\left( {t - a} \right)}^{k + n_i- \alpha_i }}} \\
=&\hspace{-6pt} \sum\nolimits_{k = 0}^{ + \infty } {\big( {\begin{smallmatrix}
\alpha_i \\
k
\end{smallmatrix}} \big)\frac{{{f^{\left( k \right)}}\left( t \right)}}{{\Gamma \left( {k - \alpha_i+ 1 } \right)}}{{\left( {t - a} \right)}^{k - \alpha_i }}} ,
\end{array}
\end{eqnarray}
for $t\in[t_{i-1},t_i)$ and
\begin{equation}\label{Eq101}
\begin{array}{l}
{}_a^{\rm{R}}{\mathscr D}_t^{\alpha\left( t \right)} f\left( t \right)= \sum\nolimits_{k = 0}^{ + \infty } {\big( {\begin{smallmatrix}
\alpha\left( t \right) \\
k
\end{smallmatrix}} \big)\frac{{{f^{\left( k \right)}}\left( t \right)}}{{\Gamma \left( {k - \alpha\left( t \right)+ 1 } \right)}}{{\left( {t - a} \right)}^{k - \alpha\left( t \right) }}} ,
\end{array}
\end{equation}
for $t\in[a,t]$.

In the case that $n-1<\alpha\left( t \right)<n$, $n\in \mathbb{N}_+$, the $\alpha\left( t \right)$-th Caputo fractional derivative of $f(t)$ can be calculated from Definition \ref{Definition 6} as
\begin{equation}\label{Eq102}
{\textstyle
\begin{array}{rl}
{}^{\rm C}_a{\mathscr D}_{t}^{\alpha\left( t \right)} f\left( t \right)=&\hspace{-6pt}{{}_a^{\rm R}}{\mathscr I}_{t}^{n-\alpha\left( t  \right)}\frac{{\rm{d}}^{{n}}}{{{\rm{d}}t^{{n}}}} f\left( t \right)\\
=&\hspace{-6pt}\sum\nolimits_{k = 0}^{ + \infty } {\frac{{{f^{\left( n+k \right)}}\left( a \right)}}{{\Gamma \left( {k +n- \alpha\left( t \right) + 1} \right)}}{{\left( {t - a} \right)}^{k +n- \alpha\left( t \right) }}}\\
=&\hspace{-6pt} \sum\nolimits_{k = n}^{+\infty}  {\frac{{{f^{\left( k \right)}}\left( a \right)}}{{\Gamma \left( { k- \alpha\left( t \right)+ 1 } \right)}}} {\left( {t - a} \right)^{k - \alpha\left( t \right) }}.
\end{array}}
\end{equation}
Additionally, one has
\begin{equation}\label{Eq103}
\begin{array}{rl}
{}_a^{\rm{C}}{\mathscr D}_t^{\alpha\left( t \right)} f\left( t \right)=&\hspace{-6pt}\frac{1}{{\Gamma \left( {n - \alpha\left( t \right) } \right)}}\int_a^t {{{\left( {t - \tau } \right)}^{n - \alpha\left( t \right)  - 1}}{f^{\left( n \right)}}\left( \tau  \right){\rm{d}}\tau } \\
=&\hspace{-6pt}\frac{1}{{\Gamma \left( {n - \alpha\left( t \right) } \right)}}\int_a^t {\sum\nolimits_{k = 0}^{ + \infty } {{{\left( { - 1} \right)}^k}\frac{{{f^{\left( {k + n} \right)}}\left( t \right)}}{{\Gamma \left( {k + 1} \right)}}} {{\left( {t - \tau } \right)}^{n + k - \alpha\left( t \right)  - 1}}{\rm{d}}\tau } \\
=&\hspace{-6pt}\frac{1}{{\Gamma \left( {n - \alpha\left( t \right) } \right)}}\sum\nolimits_{k = 0}^{ + \infty } {{{\left( { - 1} \right)}^k}} \int_a^t {\frac{{{f^{\left( {k + n} \right)}}\left( t \right)}}{{\Gamma \left( {k + 1} \right)}}{{\left( {t - \tau } \right)}^{n + k - \alpha\left( t \right)  - 1}}{\rm{d}}\tau } \\
=&\hspace{-6pt}\frac{1}{{\Gamma \left( {n - \alpha\left( t \right) } \right)}}\sum\nolimits_{k = 0}^{ + \infty } {{{\left( { - 1} \right)}^k}\frac{{{f^{\left( {k + n} \right)}}\left( t \right)}}{{\Gamma \left( {k + 1} \right)}}} \frac{{{{\left( {t - a} \right)}^{n + k - \alpha\left( t \right) }}}}{{n + k - \alpha\left( t \right) }}\\
=&\hspace{-6pt} \sum\nolimits_{k = n}^{ + \infty } {\frac{{{{\left( { - 1} \right)}^k}{f^{\left( k \right)}}\left( t \right)}}{{\Gamma \left( {n - \alpha\left( t \right) } \right)\Gamma \left( {k - n - 1} \right)\left( {k - \alpha\left( t \right) } \right)}}{{\left( {t - a} \right)}^{k - \alpha\left( t \right) }}} \\
=&\hspace{-6pt} \sum\nolimits_{k = n}^{ + \infty } {\big( {\begin{smallmatrix}
{\alpha\left( t \right)  - n}\\
{k - n}
\end{smallmatrix}} \big)\frac{{{f^{\left( k \right)}}\left( t \right)}}{{\Gamma \left( {k - \alpha(t)+ 1} \right)}}{{\left( {t - a} \right)}^{k - \alpha\left( t \right) }}}.
\end{array}
\end{equation}

\begin{remark}\label{Remark 3}
It should be noted that the variable order case is more complicated than the constant order case. To obtain concise results, some assumptions are made. Actually, more related work is expected to be done on this aspect. Additionally, the elaborated left-handed observations and discussions could adapt to the right-handed case, i.e., ${{}_t^{\rm R}}{{\mathscr I}_{b}^\alpha}{f\left( t \right)}$, ${{}_t^{\rm R}}{{\mathscr D}_{b}^\alpha}{f\left( t \right)}$, ${{}_t^{\rm C}}{{\mathscr D}_{b}^\alpha}{f\left( t \right)}$, ${{}_t^{\rm R}}{{\mathscr I}_{b}^{\alpha\left( t \right)}}{f\left( t \right)}$, ${{}_t^{\rm R}}{{\mathscr D}_{b}^{\alpha\left( t \right)}}{f\left( t \right)}$ and ${{}_t^{\rm C}}{{\mathscr D}_{b}^{\alpha\left( t \right)}}{f\left( t \right)}$. In pursuit of the clarity and brevity, the related context is not provided here.
\end{remark}

\begin{remark}\label{Remark 4}
Before ending the main results, the main contributions are summarized as follows. i) Two kinds of Taylor series like representations are proposed for Caputo fractional derivative. ii) Two kinds of Taylor series like representations are summarized for Riemann--Liouville fractional integrals/derivatives. iii) Two kinds of Taylor formula like representations are presented for fractional integrals/derivatives. iv) Some interesting properties are confirmed or newly deduced with the help of the developed representations. v) Taylor series like representations are extended from the constant order case to the variable order case. vi) New numerical approximation algorithms for fractional derivative computing are introduced carefully.
\end{remark}
\section{Simulation Study}\label{Section 4}
\begin{example}\label{Example 1}
Suppose that $f\left( t \right)$ is a non-constant periodic function with a period $T$, i.e.,
\begin{equation}\label{Eq104}
f\left( t \right)=f\left( t +T \right)
\end{equation}
for all $t\ge a$. It follows from (\ref{Eq3}) that
\begin{equation}\label{Eq105}
\begin{array}{rl}
{}_a^{\rm{C}}{\mathscr D}_{t + T}^\alpha f\left( {t + T} \right)=&\hspace{-6pt} \sum\nolimits_{k = n}^{ + \infty } {\left( {\begin{smallmatrix}
\alpha-n \\
k-n
\end{smallmatrix}} \right)\frac{{{f^{\left( k \right)}}\left( {t + T} \right)}}{{\Gamma \left( {k - \alpha+ 1 } \right)}}{{\left( {t + T - a} \right)}^{k - \alpha }}} \\
 =&\hspace{-6pt} \sum\nolimits_{k = n}^{ + \infty } {\left( {\begin{smallmatrix}
\alpha-n \\
k-n
\end{smallmatrix}} \right)\frac{{{f^{\left( k \right)}}\left( t \right)}}{{\Gamma \left( {k - \alpha+ 1 } \right)}}{{\left( {t + T - a} \right)}^{k - \alpha }}} \\
 \ne&\hspace{-6pt} \sum\nolimits_{k = n}^{ + \infty } {\left( {\begin{smallmatrix}
\alpha-n \\
k-n
\end{smallmatrix}} \right)\frac{{{f^{\left( k \right)}}\left( t \right)}}{{\Gamma \left( {k- \alpha + 1 } \right)}}{{\left( {t - a} \right)}^{k - \alpha }}} \\
 =&\hspace{-6pt} {}_a^{\rm{C}}{\mathscr D}_t^\alpha f\left( t \right),
\end{array}
\end{equation}
which indicates that for any $\alpha\in(n-1,n)$, $n\in\mathbb{N}_+$, ${}_a^{\rm{C}}{\mathscr D}_t^\alpha f\left( t \right)$ is not a periodic function anymore (see \cite{Tavazoei:2010Auto}). We here provide a way to construct a periodic fractional derivative for $f(t)$. Let $t-a\to+ \infty$, one has
\begin{equation}\label{Eq106}
\begin{array}{rl}
{\left( {t + T - a} \right)^{k - \alpha }} =&\hspace{-6pt} {\left( {t - a} \right)^{k - \alpha }}\big[ {1 + \left( {k - \alpha } \right)\frac{{T}}{{t - a}} + O\big( {\frac{{{T^2}}}{{{{( {t - a} )}^2}}}} \big)} \big]\\
 \approx&\hspace{-6pt} {\left( {t - a} \right)^{k - \alpha }} + \left( {k - \alpha } \right)T{\left( {t - a} \right)^{k - 1 - \alpha }},
\end{array}
\end{equation}
which leads to
\begin{equation}\label{Eq107}
\begin{array}{rl}
{}_a^{\rm{C}}{\mathscr D}_{t + T}^\alpha f\left( {t + T} \right) \approx&\hspace{-6pt} {}_a^{\rm{C}}{\mathscr D}_t^\alpha f\left( t \right) + \frac{T}{{t - a}}\sum\nolimits_{k = n+1}^{ + \infty } {\left( {\begin{smallmatrix}
\alpha-n \\
k-n
\end{smallmatrix}} \right)\frac{{\left( {k - \alpha } \right){f^{\left( k \right)}}\left( t \right)}}{{\Gamma \left( {k - \alpha+ 1 } \right)}}{{\left( {t - a} \right)}^{k - \alpha }}} \\
 =&\hspace{-6pt} {}_a^{\rm{C}}{\mathscr D}_t^\alpha f\left( t \right) - \frac{{T\Gamma \left( {\alpha-n  + 1} \right)\sin \left( {\alpha \pi } \right)}}{{\pi {{\left( {t - a} \right)}^{\alpha  + 1-n}}}}\sum\nolimits_{k = 1}^{ + \infty } {\frac{{{f^{\left( k +n \right)}}\left( t \right)}}{{k!}}{{\left( {a - t} \right)}^{k }}} \\
 =&\hspace{-6pt} {}_a^{\rm{C}}{\mathscr D}_t^\alpha f\left( t \right) - \frac{{T\Gamma \left( {\alpha  + 1} \right)\sin \left( {\alpha \pi } \right)}}{{\pi {{\left( {t - a} \right)}^{\alpha  + 1-n}}}}[f^{(n)}\left( a \right)-f^{(n)}\left( t \right)]\\
 \approx&\hspace{-6pt} {}_a^{\rm{C}}{\mathscr D}_t^\alpha f\left( t \right).
\end{array}
\end{equation}

It is ready to conclude that with the mentioned hypothesis, i.e., $t-a\to+\infty$, ${}_a^{\rm{C}}{\mathscr D}_t^\alpha f\left( t \right)$ is a periodic function with the same period $T$, which confirms (\ref{Eq82}) and (\ref{Eq83}).

To eliminate the effect caused by the remote historical data, fixed memory length is considered, hence
\begin{equation}\label{Eq108}
{}_{t + T - L}^{\hspace{16pt}\rm{C}}{\mathscr D}_{t + T}^\alpha f\left( {t + T} \right) = {}_{t - L}^{\hspace{7pt}\rm{C}}{\mathscr D}_t^\alpha f\left( t \right)
\end{equation}
holds for all $t>L+a$. 
\end{example}
\begin{example}\label{Example 2}
Consider a fractional order nonlinear system described by
\begin{equation}\label{Eq109}
{}_a^{\rm{C}}{\mathscr D}_t^\alpha x\left( t \right) = f\left( x\left( t \right) \right),
\end{equation}
where $x(t)\in\mathbb{R}$, $\alpha\in(0,1)$ and $f\left( \cdot\right)$ is continuous. It has been shown in \cite{Shen:2014Auto} that the equilibrium of system (\ref{Eq109}) can never be finite-time stable. We here provide a different way to obtain this conclusion, by using the series representations developed in this paper. By reduction to absurdity, suppose that there exists a finite instant $T>a$, such that $x(t)=0,~\forall t\ge T$. Then one can calculate that
\begin{equation}\label{Eq110}
\begin{array}{rl}
x\left( t \right)=&\hspace{-6pt} {}_a^{\rm R}{\mathscr I}_t^\alpha f\left( x\left( t \right) \right)+x\left( a \right)\\
 = &\hspace{-6pt}\sum\nolimits_{k = 0}^{ + \infty } {\frac{{{f^{\left( k \right)}}\left( {x\left( a \right)} \right)\dot x\left( a \right)}}{{\Gamma \left( {k + \alpha+ 1 } \right)}}{{\left( {t - a} \right)}^{k+ \alpha }}}+x\left( a \right) \\
 = &\hspace{-6pt}0.
\end{array}
\end{equation}

Eliminating the nonzero common factor, one arrives at
\begin{equation}\label{Eq111}
{\textstyle \sum\nolimits_{k = 0}^{ + \infty } {\frac{{{f^{\left( k \right)}}\left( {x\left( a \right)} \right)}}{{\Gamma \left( {k + \alpha+ 1 } \right)}}{{\left( {t - a} \right)}^{k+\alpha}}}  =- \frac{x(a)}{\dot x(a)}.}
\end{equation}

Select ${t_k} > T - a$ incrementally with the increase of $k\in \mathbb{N}_+$ and define ${z_k} = \frac{{{f^{\left( k \right)}}\left( {x\left( a \right)} \right)}}{{\Gamma \left( {k + 1 + \alpha } \right)}}$, $c=-\frac{x(a)}{\dot x(a)}$. An infinite-dimensional equation
\begin{equation}\label{Eq112}
\left[ {\begin{array}{*{20}{c}}
1&{{t_1^\alpha}}& \cdots &{t_1^{\alpha+ k - 1 }}& \cdots \\
1&{{t_2^\alpha}}& \cdots &{t_2^{\alpha+ k - 1}}& \cdots \\
 \vdots & \vdots & \ddots & \vdots & \vdots \\
1&{{t_k^\alpha}}& \cdots &{t_k^{\alpha+ k - 1}}& \cdots \\
 \vdots & \vdots & \cdots & \vdots & \ddots
\end{array}} \right]\left[ {\begin{array}{*{20}{c}}
{{z_1}}\\
{{z_2}}\\
 \vdots \\
{{z_k}}\\
 \vdots
\end{array}} \right] = \left[ {\begin{array}{*{20}{c}}
c\\
c\\
 \vdots \\
c\\
 \vdots
\end{array}} \right]
\end{equation}
follows. Observe that the coefficient matrix $A$ in (\ref{Eq112}) is a Vandermonde like matrix with full rank. As a result, (\ref{Eq112}) has a unique solution. Since the right hand of (\ref{Eq112}) is constant and $A$ is time varying, the solution of (\ref{Eq112}) must be time varying, which leads to a contradiction regarding the fact that $z_k$ is time invariant. This conclusion coincides with Theorem 7 in \cite{Shen:2014Auto}.
\end{example}

\begin{example}\label{Example 3}
To evaluate the approximation provided so far, the exponential function ${\rm e}^{2t}$ is considered as a test function. After computing the analytic solution via (\ref{Eq76}) with $\alpha=0.5$ and $a=0$, the comparison results are illustrated in Fig. 1 and Fig. 2, where $E$ is introduced for quantitative analysis defined as the norm of approximation error and the sampling period is $0.001$ second.

The two figures show that i) both (\ref{Eq84}) and (\ref{Eq86}) are effective approximations of Riemann--Liouville fractional derivative; ii) a larger $N$ leads to a smaller approximation error; and iii) compared with (\ref{Eq84}), (\ref{Eq86}) performs better with the same $N$.

Next, we will check the approximation performance on Caputo fractional derivative. With the relationship in (\ref{Eq9}), (\ref{Eq85}) is implied in (\ref{Eq84}) and Fig. 1. Consequently, only (\ref{Eq87}) is studied with the related results shown in Fig.  3. It can be observed that such an approximation gives a satisfied degree of accuracy to that of (\ref{Eq84}), while it is slightly poorer than using (\ref{Eq86}). Another interesting point should be emphasized is that the Caputo fractional derivative of ${\rm e}^{2t}$ is monotonic while the Riemann--Liouville one is not. The reason is that ${}_0^{\rm{R}}{\mathscr D}_t^{0.5} {\rm e}^{2t}$ includes the decreasing component $\frac{2^k}{\Gamma(k+0.5)}t^{k-0.5}$ with $k<n$ and the increasing component $\frac{2^k}{\Gamma(k+0.5)}t^{k-0.5}$ with $k>n$ simultaneously.
\end{example}

\section{Conclusions}\label{Section 5}
In this paper, two kinds of infinite series composed of integer derivatives and power functions have been investigated for fractional derivative and fractional integral systematically. With the series representation established, several essential properties are established, such as singularity, additivity and differentiability, etc. Besides, theoretical results are vividly validated by numerical examples. It is anticipated that the elaborated results would enrich the understanding and utilizing of fractional calculus. The following interesting and promising topics should be investigated in the foreseeable future.
\begin{enumerate}[i)]
  \item Apply the truncated series for numerical computation of fractional differential equations.
  \item Extend the results to the fractional integrals/derivatives of a function with distributed order.
  \item Study the series representation of fractional derivatives for functions that are not analytic.
\end{enumerate}

\vspace*{-7pt}

\section*{Acknowledgements}
The authors would like to thank the Associate Editor and the anonymous reviewers for their keen and insightful comments which greatly improved the contents and the presentation of this paper. The work in this paper was supported by the National Natural Science Foundation of China (61601431, 61573332), the Anhui Provincial Natural Science Foundation (1708085QF141) and the fund of China Scholarship Council (201806345002).

\section*{References}
\bibliographystyle{model1-num-names}
\bibliography{database}


\end{document}